\documentclass[12pt,journal]{new-aiaa}
\usepackage[utf8]{inputenc}
\usepackage{textcomp}
\usepackage{ulem} 
\usepackage{graphicx}
\usepackage{amsmath}
\usepackage[version=4]{mhchem}
\usepackage{siunitx}
\usepackage{longtable,tabularx}
\usepackage{color}
\usepackage{mathtools}
\usepackage{float}
\usepackage{caption,subcaption}
\usepackage{tikz}
\usepackage[many]{tcolorbox}
\usepackage{color}
\usepackage{epstopdf}
\usepackage[numbers,sort&compress]{natbib}
\usepackage{algorithm}
\usepackage{algorithmic}
\usepackage{setspace}
\usepackage[margin=1in]{geometry}
\def\x{\boldsymbol{x}}

\def\1{\boldsymbol{1}}

\newtheorem{problem}{Problem}

\title{Fuel-optimal powered descent guidance for lunar pinpoint landing using neural networks}%

\author[1]{Kun Wang\footnote{PhD student, School of Aeronautics and Astronautics, Student member.}}
\author[1,2*]{Zheng Chen\footnote{Researcher, School of Aeronautics and Astronautics. email: \underline{z-chen@zju.edu.cn} (Corresponding author). }}
\author[1,2]{Jun Li\footnote{Professor, School of Aeronautics and Astronautics.}}
\affil[1]{School of Aeronautics and Astronautics, Zhejiang University, Hangzhou 310027, Zhejiang, China}
\affil[2]{Huanjiang Lab, Zhuji, 311816, Zhejiang, China}

\begin{document}
\maketitle
\begin{abstract}
This paper presents a Neural Networks (NNs) based approach for designing the Fuel-Optimal Powered Descent Guidance (FOPDG) for  lunar pinpoint landing. According to the Pontryagin's Minimum Principle, the optimality conditions are firstly derived. To generate the dataset of optimal trajectories for training NNs, we formulate a parameterized system, which allows for generating each optimal trajectory by a simple propagation without using any optimization method. Then, a dataset containing the optimal state and optimal thrust vector pairs can be readily collected. Since it is challenging for NNs to approximate bang-bang (or discontinuous) type of optimal thrust magnitude, we introduce a regularisation function to the switching function so that the regularized switching function approximated by a simple NN can be used to represent the optimal thrust magnitude.
Meanwhile, another two well-trained NNs are used to predict the thrust steering angle and time of flight given a flight state. Finally, numerical simulations show that the proposed method is capable of generating the FOPDG that steers the lunar lander to the desired landing site with acceptable landing errors.
\end{abstract}




\section{Introduction}
\label{intro}
As the world is becoming increasingly interested in lunar exploration, it is expected that each lunar lander can automatically steer itself to a predefined site during the powered descent phase. The guidance command, which determines the control (or thrust vector), plays a key role in guaranteeing an automatic and safe landing. Hence, significant efforts have been devoted in both academia and industry to the Powered Descent Guidance (PDG), which can be roughly categorized into three classes: Analytical PDG (APDG), optimization based PDG, and learning based PDG \citep{song2020survey}. 


Since the computational ability was limited in the Apollo era, 
it is necessary to devise APDG by simplifying the pinpoint landing problem. 
 The well-known Apollo PDG was developed by representing the thrust vector as a linear function of time \citep{cherry1964general} (also known as E-guidance), or a quadratic polynomial of time \citep{klumpp1974apollo}, depending on different terminal constraints. As a result, the analytical expressions for speed and position could be obtained on integration of dynamics equation. The coefficients in the polynomials were then determined based on the boundary conditions at the predefined landing site. The Apollo PDG is well known to require few efforts to implement as it has a closed-form solution. Other works have also been proposed based on Apollo PDG \citep{lu2019augmented,lu2020theory}. On the other hand, the Apollo PDG has some drawbacks, including the need to estimate the time of flight, loss of optimality in terms of fuel consumption, and inability to handle limits on thrust magnitude. In addition, some zero-effort-miss/zero-effort-velocity based feedback guidance algorithms have been proposed recently for the Mars pinpoint landing problem
 \citep{zhou2014improved,zhang2016multi,zhang2017collision,wang2021two}. One common assumption made in those works is that the gravitational field is uniform. However, the flight time and flight range for the powered descent at the Moon are usually longer than that at Mars. Therefore, the assumption of uniform gravitational field is not applicable to the lunar pinpoint landing. In this case, the gravity variation during the powered descent is highly nonlinear, which is usually captured by the standard Newtonian gravity model. This will become a nuisance when devising the analytical PDG. 


Due to limited capacity of the lunar lander, devising the Fuel-Optimal PDG (FOPDG) for the pinpoint landing problem has been attracting plenty of attention in recent decades. The FOPDG 
can be formulated as an Optimal Control Problem (OCP), which is usually solved via indirect or direct methods \citep{rao2009survey}. The indirect methods are usually based on the Pontryagin's Minimum Principle (PMP) to convert the OCP to a Two-Point Boundary-Value Problem (TPBVP). Then, shooting methods can be used to solve the resulting TPBVP. In this regard, 
\citep{topcu2007minimum} investigated the properties of the FOPDG and found that the optimal thrust magnitude is  bang-bang for both point-mass and rigid-body lander models. In addition, the optimal thrust magnitude was demonstrated to have at most two switches for a two-dimensional powered descent problem. For a three-dimensional FOPDG problem, 
\citep{lu2018propellant} modeled the thrust magnitude switching structure as a maximum–minimum–maximum profile, where the first switching time was determined empirically and the second switching time was determined through a univariant optimization. However, 
determining the switching times by the use of switching function could be difficult, because it needs the mass costate for calculating the switching function. To this end, 
\citep{ito2020throttled} proposed a simple approximation to the switching function without the need for the mass costate. Since the gradients of switching times have discontinuities, resulting in potential numerical difficulties, a sigmoid function was employed to smooth the bang-bang thrust magnitude \citep{lu2023propellant}. The homotopy technique was embedded into the indirect method to find the fuel-optimal trajectory for landing on an asteroid with a highly irregular gravitational field
\citep{yang2020fuel}. In addition, some works on devising the FOPDG under complex constraints, such as glide slope constraint \citep{you2022theoretical} and thrust pointing constraint \citep{ito2023optimal,lu2023propellant} are also noteworthy. 

 
Although indirect methods are widely used to solve OCPs, providing an appropriate initial guess for the costate vector is intricate. Additionally, the optimal control profile should be given {\it a priori} in the presence of state and/or control constraints. In contrast, direct methods transform the OCP into a Non-Linear Programming (NLP) problem via collocation methods, which is then solved by interior-point or sequential quadratic programming method \citep{rao2009survey}. Because of its efficiency in dealing with path constraints, direct methods have been widely employed to solve OCPs. However, the resulting optimization problem can become extremely challenging to solve. Thanks to the recent development of numerical programming and hardware, convex optimization has been becoming a powerful approach for solving OCPs in real time regarding the soft landing \citep{accikmecse2013lossless,liu2019fuel,SaglianoRocketLanding2023}. In \citep{acikmese2007convex}, a 3 Degree-of-Freedom (DoF) PDG problem with nonconvexity resulting from the lower bound on the thrust magnitude was transformed to a convex second-order cone programming problem. The resulting problem could be solved in polynomial time. Later, this lossless convexification methodology was extended to PDG problems with more complex constraints, such as  thrust pointing constraint \citep{accikmecse2013lossless} and state constraints \citep{dueri2017customized}. Moreover, a real-time successive convexification algorithm for a generalized free-final-time 6-DoF PDG problem with state-triggered constraints was presented \citep{szmuk2020successive}. The advantage of convex-optimization-based FODPGs is that they can generate a fuel-optimal trajectory subject to complex constraints. However, a tailored convexification method may be required for a specified problem, and a customized solver may be in demand for real-time implementation \citep{elango2022customized}. In addition, by transforming a dynamic programming problem  into a static optimization problem, an approach called model predictive static programming has been proposed to devise the FODPG \citep{sachan2015fuel}. Nevertheless, this approach may fail to converge when the initial guess is not sufficiently close to the solution \citep{pan2019newton}.

The techniques mentioned in the last two paragraphs are widely used, but they are
usually time-consuming and may suffer from convergence issues \citep{chen2019nonlinear}. 
Therefore, the above methods are typically not suitable for onboard systems with limited computational resources. Thanks to the emerging achievements of artificial intelligence in recent decades,  using machine learning to generate real-time solutions for OCPs has been attracting great attention for aerospace applications, ranging from orbital transfer \citep{izzo2021real,wang2023real}, soft landing \citep{sanchez2018real,cheng2019real,gaudet2020deep,you2021learning,song2021feasibility,zhao2022real,wangJinbo2023real}, to missile guidance \citep{wang2022nonlinear,dai2023entry}.  These works are mainly based on two frameworks, i.e., reinforcement learning and Supervised Learning (SL). Based on the universal approximation theorem \citep{HORNIK1989359}, Neural Networks (NNs), especially deep NNs, have been widely used for approximating highly nonlinear mappings within both frameworks. Regarding the existing works based on the SL, it usually requires to solve a large number of OCPs with different boundary conditions via either indirect \citep{sanchez2018real,cheng2019real} or direct methods \citep{you2021learning}. Then, the optimal state-control or optimal state-costate pairs are stored in a dataset, which is necessary for training NNs. However, both indirect and direct methods may fail to converge, so the process of dataset generation could be tedious. In addition, the universal approximation theorem based NNs are only able to approximate a Borel measurable function, i.e., a continuous function \citep{HORNIK1989359}. Recall that the optimal thrust magnitude is bang-bang for the FOPDG problem, so oscillations and under- or over-shooting are expected to appear in the thrust magnitude prediction generated from well-trained NNs, as shown by the numerical simulations in \citep{cheng2019real,wangJinbo2023real}. In this case, approximating the optimal thrust magnitude directly can be quite challenging \citep{DeeperLearning2024}, and 
the resulting thrust magnitude may be unimplementable, which could cause critical issues for the pinpoint landing mission.

In this paper an SL-based method to generate the optimal thrust vector in real time for the lunar pinpoint landing is developed.   To generate the training dataset,
instead of using  conventional optimization-methods that suffer from convergence issues, we formulate a parameterized system based on the necessary conditions in virtue of PMP \citep{wang2022nonlinear,wang2023real}. This allows for generating an optimal trajectory by a simple propagation on the parameterized system. In order to avoid the challenging issue of using NNs to approximate the  bang-bang thrust magnitude, we attempt to use NNs to approximate the continuous switching function instead. Nevertheless, we show that the relationship between the flight state and the switching function exhibits a feature of one-to-many mapping, containing contradictory input-output pairs \citep{izzo2020stability}. In such case, NNs are not applicable to approximating such mapping since they are only able to learn a one-to-one mapping \citep{li2022using}. To resolve this issue, a regularisation function is introduced into the switching function so that the one-to-many mapping can be eliminated. As a result, the approximation performance of NNs is greatly enhanced, and the regularised switching function can be approximated precisely. Thus, the thrust magnitude can be derived accordingly given a flight state. Meanwhile, another two NNs are built and trained, allowing for generating the thrust steering angle and predicting the time of flight.

The remainder of the paper is organized as follows. The OCP for the lunar pinpoint landing is formulated  in Section \ref{SE:problem}. In Section \ref{SE:properties}, optimality conditions are derived according to PMP, and a system is formulated to parameterize the optimal trajectory. Section \ref{SE:Real} describes the scheme for generating the FOPDG, including the dataset generation, regularisation procedure for the switching function, and the training process of NNs. Section \ref{SE:Numerical} presents some numerical examples, demonstrating and verifying the developments of the paper. This paper finally concludes by Section \ref{SE:conclusions}. 


\section{Problem Formulation}\label{SE:problem}
\begin{figure}[!htp]
\begin{center}
\includegraphics[scale=0.28]{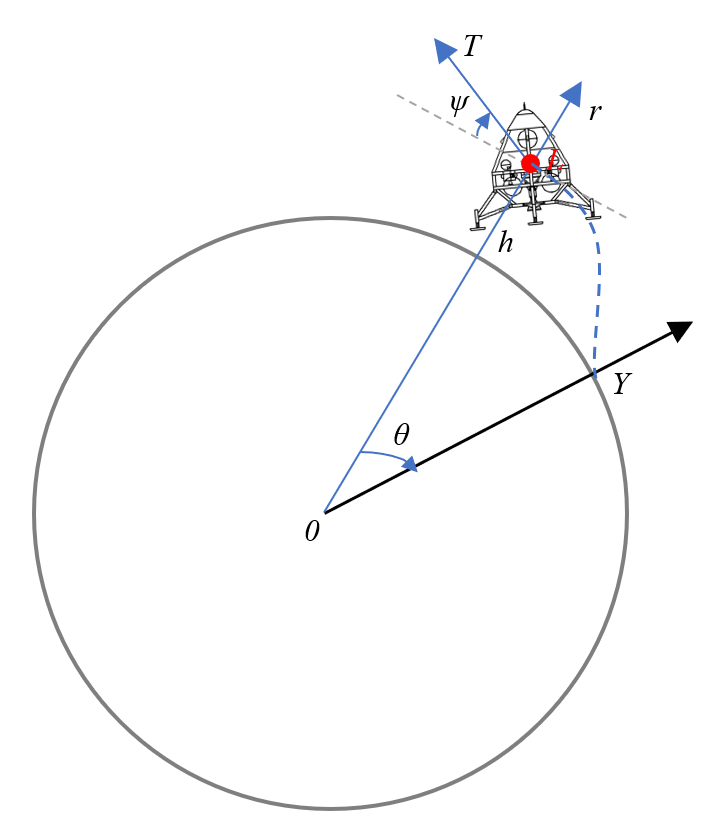}
\caption{Coordinate system for the lunar pinpoint landing.}\label{Fig:frame}
\end{center}
\end{figure}

Consider the planar motion of a lunar lander, as presented in Fig.~\ref{Fig:frame}. 
It is assumed that the Moon is a regular spherical body and we ignore the influences of the Moon’s rotation. The point $O$ is located at the center of the Moon, and the desired landing site on the lunar surface is specified by the point $Y$.
Let $r \ge R_0$ ($R_0$ is the radius of the Moon) denote the radial distance between the point $O$ and the lunar lander $L$, thus the altitude of the lunar lander is  $h = r - R_0$. The angle between the line segment $OL$ and $OY$ is called range angle, denoted by $\theta\in [0,2\pi]$. The lunar lander is propelled by an engine with adjustable magnitude and steering angle. Denote by $T \in [0,T_{m}]$the thrust magnitude, and the constant $T_{m}$ is the maximum thrust magnitude. Let $u\in[0,1]$ denote the engine thrust ratio, then we have $T=uT_{m}$. The thrust steering angle $\psi\in [-\pi, \pi]$ is defined to be the angle from the local horizontal line of the lunar lander to the thrust vector. Then, the point-mass dynamics of the lunar lander is given by \citep{liu2008optimal}
\begin{align}
\begin{cases}
\dot{r}(t) =  v(t),\\
\dot{v}(t) =  \frac{u(t)T_{m}}{m(t)}\sin\psi(t)-\frac{\mu}{r^2(t)} + r(t)\omega^2(t),\\
\dot{\theta}(t) = -\omega(t), \\
\dot{\omega}(t) = -[\frac{u(t)T_{m}}{m(t)}\cos\psi(t) + 2v(t)\omega(t)]/r(t),\\
\dot{m}(t) = -\frac{u(t)T_{m}}{I_{sp}g_e},
\end{cases}
\label{LunarLander:DiffEqution}
\end{align}
where $t \ge 0$ is the time, $v$ the radial speed along the direction of $OL$, and $\omega$  the angular velocity of the lunar lander. Thus, the transverse speed is $\omega r$. $\mu$ is the gravitational constant of the Moon. $m$ is the mass of the lunar lander. The constant $I_{sp}$ denotes the specific impulse of the lunar lander's propeller, and $g_e$ represents the Earth's gravitational acceleration at sea level. Let $t_0 = 0$ be the initial time, and the 
initial condition of the lunar lander is fixed, i.e.,
\begin{align}
r(0) = r_0, ~v(0) = v_0, ~\theta(0) = \theta_0, ~\omega(0) = \omega_0, ~m(0) = m_0.
\label{InitialState}
\end{align}

From Fig.~\ref{Fig:frame}, at touchdown we have $r(t_f) = R_0$ ($t_f$ is the free final time), and the landing site for the lunar lander can be characterized by the final state of the range angle, i.e., $\theta_f$. For a successful pinpoint landing, we expect to have zero speed at touchdown, namely $v(t_f) = 0$ and $\omega(t_f) = 0$. Then, guiding the lunar lander to the predefined landing site, while consuming minimum fuel, requires addressing the following OCP in real time.
\begin{problem}\label{problem:OIP}
Given an initial condition defined in Eq.~(\ref{InitialState}) and a predefined landing site, solving the OCP is equivalent to finding the guidance command $u(t)$ and $\psi(t)$,
to steer the system in Eq.~(\ref{LunarLander:DiffEqution}) from the initial condition to the final condition given by 
\begin{align}
r(t_f) = R_0, ~v(t_f) = 0, ~\theta(t_f) = 0, ~\omega(t_f) = 0,
\label{FinalState}
\end{align}
such that the cost functional
$ J = \int_0^{t_{f}} u(t)~\mathrm{d}t$
is minimized. 
\end{problem}

It is clear that the cost functional amounts to $J = -m(t_f)$, 
which is equivalent to minimizing the fuel consumption. 
To improve the numerical conditioning, 
we use $R_0$, 
$\sqrt{\frac{\mu}{R_0}}$, 
$m_0$, $\sqrt{\frac{{R_0}^3}{\mu}}$, 
and $\frac{m_0 \mu}{{R_0}^2}$ to normalize variables $r$, $v$, $m$, $t$, 
and $T_m$ in Eq.~(\ref{LunarLander:DiffEqution}), respectively. 
As a result, the constant $\mu$ is normalized to $1$. To avoid abuse of notation, we still use the same notation as in Eq.~(\ref{LunarLander:DiffEqution}) for the dimensionless counterpart in the remainder of the paper.

Although conventional optimization-based methods are available for solving Problem \ref{problem:OIP}, they are not suited for onboard implementation. In the following two sections, we shall show how to develop a real-time method for generating the FOPDG by combining PMP with NNs.

\section{Parametrization of Optimal Trajectories}\label{SE:properties}

In this section, the necessary conditions for optimality will be first established. Then, a parameterized system is formulated so that a large number of optimal trajectories could be readily obtained, which will be used for training NNs.

\subsection{Necessary Conditions for Optimality} 
Denote by $\boldsymbol{x} = [r,v,\theta,\omega,m]^\top$ the state vector. Then, the nonlinear dynamics in Eq.~(\ref{LunarLander:DiffEqution}) can be rewritten as
\begin{align}
\dot{\boldsymbol x}(t) = \boldsymbol f(\boldsymbol x, u, \psi, t),
\label{EQ:rewrite}
\end{align}
where $\boldsymbol f:\mathbb{R}^5 \times \mathbb{R} \times \mathbb{R} \times \mathbb{R}_0^+ \rightarrow \mathbb{R}^5$ is the smooth vector field defined in Eq.~(\ref{LunarLander:DiffEqution}).

Let $\boldsymbol{p_x} = [p_r, p_v, p_{\theta}, p_{\omega},p_m]^\top $ denote the co-state vector related to $\boldsymbol{x}$. Then, 
the Hamiltonian is expressed as
\begin{align}
\mathscr{H}= {\boldsymbol f}^T \boldsymbol{p_x} + u = 
p_r v + p_v (\frac{uT_{m}}{m}\sin\psi-\frac{\mu}{r^2} + r\omega^2) - p_{\theta} \omega + p_\omega[-(\frac{uT_{m}}{m}\cos\psi + 2v\omega)/r] + p_m(-\frac{uT_{m}}{I_{sp}g_e}) + u. 
\label{EQ:Ham}
\end{align}
According to PMP \citep{Pontryagin}, along an extremal trajectory, we have    
\begin{align}
\begin{cases}
\dot{p}_r(t) = -\frac{\partial \mathscr{H}}{\partial r} =  -\frac{2p_v\mu}{r^3} -  p_v{\omega}^2 -p_{\omega}[(\frac{uT_{m}}{m}\cos\psi + 2v\omega)/r^2],\\
\dot{p}_v(t) = -\frac{\partial \mathscr{H}}{\partial v} =   -p_r + \frac{2p_{\omega}\omega}{r},\\
\dot{p}_{\theta}(t) = -\frac{\partial \mathscr{H}}{\partial \theta} = 0,\\
\dot{p}_{\omega}(t) = -\frac{\partial \mathscr{H}}{\partial \omega} = -2 p_v r\omega + p_{\theta} + \frac{2 p_{\omega} v}{r}, \\
\dot{p}_m(t) = -\frac{\partial \mathscr{H}}{\partial m} = \frac{p_v u T_{m}}{m^2} \sin\psi - \frac{p_{\omega} u T_{m}}{m^2 r} \cos\psi.
\end{cases}
\label{EQ:dot_p}
\end{align}
For the optimal thrust steering angle, it follows that
\begin{align}
\frac{\partial \mathscr{H}}{\partial \psi} = 0. \label{EQ:direction}
\end{align}
Explicitly rewriting Eq.~(\ref{EQ:direction}) leads to
\begin{align}
\psi(t) =  \arctan [-\frac{p_v(t)r(t)}{p_\omega(t) }].
\label{EQ:optimal_psi}
\end{align}
Minimizing $\mathscr{H}$ w.r.t. $\psi$ further implies that
\begin{align}
\left[ \begin{array}{c}
\sin\psi(t) \\
\cos\psi(t)
\end{array}
\right]  = - \frac{1}{\sqrt{{p^2_v(t)}+[-\frac{p_w(t)}{r(t)}]^2}}\left[ \begin{array}{c}
p_v(t) \\
-\frac{p_w(t)}{r(t)}
\end{array}
\right].
\label{EQ:optimal_H1}
\end{align}
We assume that the optimal engine thrust ratio is bang-bang, i.e.,
\begin{align}
u(t)= 
\left\{ 
    \begin{array}{lc}
        1, S(t) < 0 \\
        0, S(t) > 0\\
    \end{array}
\right.
\label{EQ:magnitudenew}
\end{align}
where $S(t)$ denotes the switching function:
 \begin{align}
S(t) = 1 - \frac{T_{m}}{m(t)} \sqrt{p^2_v(t) + [\frac{p_\omega(t)}{r(t)}]^2 }-\frac{p_m(t)T_{m}}{I_{sp}g_e}.
\label{EQ:SF}
\end{align}
As the final mass of the lunar lander is free, the transversality condition implies
\begin{align}
p_{m}(t_f) = 0.
\label{EQ:Transversality1}
\end{align}
Because the Hamiltonian in Eq.~(\ref{EQ:Ham}) does not contain time explicitly and the final time is free, we have 
\begin{align}
\mathscr{H}(t)\equiv 0,\ \forall~ t\in[0,t_f].
\label{EQ:hamiszero}
\end{align}
For brevity, a quintuple $(r(t), v(t), \theta(t), \omega(t), m(t))$ for $t \in [0, t_f]$ is said to
be an optimal trajectory if all the necessary conditions in Eqs.~(\ref{EQ:dot_p}),(\ref{EQ:optimal_psi}),(\ref{EQ:magnitudenew}), (\ref{EQ:Transversality1}) and (\ref{EQ:hamiszero}) are
met. In order to generate the FOPDG in real time via NNs, a training dataset containing a large number of optimal trajectories is required. In this regard, one
viable approach is to employ some root-finding algorithms to solve the following shooting function:
\begin{align}
\boldsymbol{\Phi} (\boldsymbol{{p_x}_0},t_f) = [r(t_f)-R_0,v(t_f),\theta(t_f), \omega(t_f),p_{m}(t_f),\mathscr{H}(t_f)]=\boldsymbol{0},
\label{EQ:TPBVP_law}
\end{align}
where $\boldsymbol{{p_x}_0}$ is the initial guess for the costate vector $\boldsymbol{p_x}$ and $t_f$ is the initial guess for the final time. However, solving Eq.~(\ref{EQ:TPBVP_law}) could be extremely challenging because $\boldsymbol{{p_x}_0}$ has no physical significance and the shooting method is very sensitive to the initial guess. Moreover, the final state for the pinpoint landing is fully constrained, making it even more difficult to solve \citep{lu2023propellant}. In the next subsection, we will present a parameterized system so that an optimal trajectory can be readily obtained.
\subsection{Parametrized System}\label{SE:Pontryagin}
Define a new independent variable $\tau$ as below
\begin{align}
\tau = t_f - t, t \in [0,t_f].
\label{Eq:tao}
\end{align}
Let us establish a first-order ordinary differential system 
\begin{align}
\begin{cases}
\dot{\bar {\boldsymbol x}} = -\boldsymbol f(\boldsymbol {\bar x}, {\bar u},  \bar \psi, \tau), \\
\dot{\bar{\boldsymbol{p}}}_{\boldsymbol x} = \frac{\partial {\bar {\mathscr{H}}}(\tau)}{\partial \boldsymbol {\bar x}(\tau)},
\label{EQ:pareEquation}
\end{cases}
\end{align}
where $\boldsymbol {\bar x} = [\bar r,\bar v,\bar \theta,\bar \omega,\bar m]^\top$, $\bar{\boldsymbol{p}}_{\boldsymbol x} =  [\bar{p}_r,\bar{p}_v,\bar {p}_\theta,\bar {p}_\omega,\bar {p}_m]^\top$,  and $\bar {\mathscr{H}}$ is defined as 
\begin{align}
\bar {\mathscr{H}} = {\bar u} + \boldsymbol{f}^\top(\boldsymbol {\bar x}, \bar u, \bar \psi)\bar{\boldsymbol{p}}_{\boldsymbol x}.
\label{Eq:H_para}
\end{align}
Meanwhile, $\bar u$ and $\bar \psi$ satisfy  
\begin{align}
\left[ \begin{array}{c}
\sin\bar{\psi}(\tau) \\
\cos\bar{\psi}(\tau)
\end{array}
\right]  = - \frac{1}{\sqrt{{\bar{p}_v}^2(\tau)+[-\frac{\bar{p}_w(\tau)}{\bar r(\tau)}]^2}}\left[ \begin{array}{c}
\bar{p}_v(\tau) \\
-\frac{\bar{p}_w(\tau)}{\bar {r}(\tau)}
\end{array}
\right],
\label{EQ:optimal_psi_para}
\end{align}
and 
\begin{align}
\bar u(\tau)= 
\left\{ 
    \begin{array}{lc}
        1, \bar {S}(\tau) < 0 \\
        0, \bar {S}(\tau) > 0\\
    \end{array}
\right.
\label{EQ:magnitudenew_system}
\end{align}
where $\bar {S}(\tau)$ is set as
\begin{align}
\bar {S}(\tau) = 1 - \frac{T_{m}}{\bar m(\tau)} \sqrt{\bar {p}^2_v(\tau) + [\frac{\bar {p}_\omega(\tau)}{\bar r(\tau)}]^2 }-\frac{\bar {p}_m(\tau)T_{m}}{I_{sp}g_e}.
\label{EQ:sf_para_system}
\end{align}
Define the initial condition at $\tau = 0$ for the system in Eq.~(\ref{EQ:pareEquation}) as
\begin{align}
\begin{cases}
\bar r(0) = 1, ~\bar v(0) = 0, ~\bar \theta (0) = 0, ~\bar \omega(0) = 0, ~\bar m(0) = \bar {m}_0, \\
\bar{p}_r(0) = \bar{p}_{r_0},~\bar{p}_v(0) = \bar{p}_{v_0},~\bar{p}_\theta(0) = \bar{p}_{\theta_0},~\bar{p}_\omega(0) = \bar{p}_{\omega_0},~\bar{p}_m(0) = 0,
\label{EQ:Initial_pareEquation}
\end{cases}
\end{align}
where the quadruple $(\bar {p}_{r_0}, \bar {p}_{v_0}, \bar {p}_{\theta_0}, \bar {p}_{\omega_0})$ is arbitrary. The value for $\bar {m}_0$ is chosen to satisfy  
\begin{align}
\bar {\mathscr{H}}(0) = 0.
\label{EQ:Hpara}
\end{align}
Substituting Eq.~(\ref{EQ:Initial_pareEquation}) into Eq.~(\ref{EQ:Hpara}) leads to
\begin{align}
\bar {\mathscr{H}}(0) = -\frac{\bar {u} T_{m}}{\bar {m}_0} \sqrt{\bar{p}^2_{v_0} + \bar {p}^2_{\omega_0}}  -\mu \bar{p}_{v_0}  + \bar {u} =0. 
\label{EQ:Ham_para_initial}
\end{align}
It is clear that the initial condition of $\boldsymbol {\bar x}$ in Eq.~(\ref{EQ:Initial_pareEquation}) can represent the normalized final condition in Eq.~(\ref{FinalState}). Therefore, it can be used to denote the final state needed for the pinpoint landing. In the final phase of the powered descent, since the engine thrust ratio is fixed at the maximum \citep{lu2023propellant}, we have $\bar {u}(0)=1$. Recall that the normalized gravitational constant $\mu$ is 1.  
Thus, Eq.~(\ref{EQ:Ham_para_initial}) further reduces to
\begin{align}
\bar {\mathscr{H}}(0) = -\frac{T_{m}}{\bar {m}_0} \sqrt{\bar{p}^2_{v_0} + \bar {p}^2_{\omega_0}}  - \bar{p}_{v_0}  + 1 =0. 
\label{EQ:Ham_para_initial_reduced}
\end{align}
Therefore, for any given pair $(\bar {p}_{v_0}, \bar {p}_{\omega_0})$, the value for $\bar {m}_0$ can be determined analytically, i.e.,
\begin{align}
\bar {m}_0 = \frac{T_{m}}{1- \bar{p}_{v_0}} \sqrt{\bar{p}^2_{v_0} + \bar {p}^2_{\omega_0}}. 
\label{EQ:Ham_para_m0}
\end{align}

Up to now, by choosing a quadruple $(\bar {p}_{r_0}, \bar {p}_{v_0}, \bar {p}_{\theta_0}, \bar {p}_{\omega_0})$ arbitrarily, the only variable $\bar {m}_0$ in Eq.~(\ref{EQ:Initial_pareEquation}) can be immediately obtained according to Eq.~(\ref{EQ:Ham_para_m0}). For the sake of notational simplicity, let the pair
\begin{align*}
(\boldsymbol {\bar x}(\tau,\bar {p}_{r_0}, \bar {p}_{v_0}, \bar {p}_{\theta_0}, \bar {p}_{\omega_0}),\bar{\boldsymbol{p}}_{\boldsymbol x}(\tau,\bar {p}_{r_0}, \bar {p}_{v_0}, \bar {p}_{\theta_0}, \bar {p}_{\omega_0})) \in \mathbb{R}^{10}
\end{align*}
for $\tau \in [0,t_f]$ be the solution to the parameterized system in Eq.~(\ref{EQ:pareEquation}) with the initial condition in Eq.~(\ref{EQ:Initial_pareEquation}).

It is clear that $\boldsymbol {\bar x}(\tau) = \boldsymbol {x}(t_f - t)$  and $\boldsymbol {\bar{p}_x}(\tau) = \boldsymbol {p_x}(t_f - t)$
hold if the quadruple $(\bar {p}_{r_0}, \bar {p}_{v_0}, \bar {p}_{\theta_0}, \bar {p}_{\omega_0})$ is chosen so that the initial condition in Eq.~(\ref{InitialState}) is the same as the final condition $\boldsymbol {\bar x}(t_f)$ obtained by propagating the parameterized system in Eq.~(\ref{EQ:pareEquation}). Meanwhile, it is palpable that the pair meets all the necessary conditions for an optimal trajectory, and $\tau$ represents the optimal time of flight. In other words, an optimal trajectory can be 
obtained by arbitrarily choosing a quadruple $(\bar {p}_{r_0}, \bar {p}_{v_0}, \bar {p}_{\theta_0}, \bar {p}_{\omega_0})$ and propagating the parameterized system in Eq.~(\ref{EQ:pareEquation}) with the initial condition in Eq.~(\ref{EQ:Initial_pareEquation}). 

Meanwhile, $\bar{\psi}$ and ${\bar u}$  can be determined through Eq.~(\ref{EQ:optimal_psi_para}) and Eq.~(\ref{EQ:magnitudenew_system}),respectively. As $\bar u$ is bang-bang, it is intricate for NNs to approximate, even with many hidden layers and/or neurons \citep{origer2023guidance}. To address this issue, we propose to approximate the continuous function ${\bar S}$ instead, which can be obtained via Eq.~(\ref{EQ:sf_para_system}). Denote by $f_\tau$ the nonlinear mapping $\boldsymbol{\bar x} \longmapsto \tau$, by $f_{\bar\psi}$ the mapping $\boldsymbol{\bar x} \longmapsto \bar \psi$, and by $f_{\bar{S}}$ the mapping $\boldsymbol{\bar x} \longmapsto \bar {S}$. 
According to the universal approximation theorem \citep{HORNIK1989359}, if we have access to a large number of sampled data representing the relationships $f_\tau$, $f_{\bar\psi}$, and $f_{\bar{S}}$, well-trained NNs will be able to accurately represent them. In the next section, we shall present the procedure for generation the FOPDG in real time by using the parameterized system and NNs.

\section{Real-Time Generation of the FOPDG via NNs}\label{SE:Real}
In this section, we begin by outlining the procedure for generating the training dataset.  Nevertheless, it is shown that the mapping $f_{\bar{S}}$ is set-valued, which prevents one from using NNs to produce an accurate and robust prediction \citep{li2022using}. To this end, a regularisation function is introduced into the switching function. Finally, three NNs are trained for real-time implementation of the FOPDG. 
\subsection{Nominal Trajectory and Dataset Generation}
Consider a lunar lander with  a dry mass of $m_{d} = 250$ kg travelling from an initial condition set as: $r_0 = 1753$ km, $v_0 = 0 $ m/s, $\theta_0 = 30 $ deg, $\omega_0 =9.6410 \times10^{-4}$ rad/s, $ m_0 =600$ kg \citep{liu2008optimal}. The propulsion system of the lunar lander is specified by $I_{sp} = 300$ s and $T_m = 1,500$ N. The radius of the Moon is $R_0 = 1738$ km, and $g_e$ is equal to 9.81 $\rm m/s^2$. Additionally, the gravitational constant $\mu$ is $4.90275 \times 10 ^{12}$~$\rm m^3/s^2$.

We solve this OCP using the indirect methods. Notice that the optimal engine thrust ratio in Eq.~(\ref{EQ:magnitudenew}) is bang-bang, which could result in numerical difficulties. In this case, some smoothing techniques \citep{bertrand2002new,wang2023new} can be used to transform the discontinuous control into a smooth one. With the smoothing technique \citep{wang2023new}, a smoothing constant $\delta = 1\times 10^{-10}$ is used to approximate the discontinuous optimal thrust engine ratio in Eq.~(\ref{EQ:magnitudenew}), i.e.,
\begin{align}
u (t)\approx u(t,\delta) = \frac{1}{2}(1-\frac{S(t)}{\sqrt{\delta+|S(t)|^2}}).
\label{EQ:approximation_u}
\end{align}
As a result, the {\it nominal trajectory} is obtained by solving Eq.~(\ref{EQ:TPBVP_law}).



Denote by $p^*_{r_f}$, $p^*_{v_f}$, $p^*_{\theta_f}$, and $p^*_{\omega_f}$ the final values of the co-states $p_{r}$, $p_{v}$, $p_{\theta}$, and $p_{\omega}$ for the {\it nominal trajectory}. Consider a new set of values given by
\begin{align}
\bar{p}_{r_0} = p^*_{r_f} + \lambda \bar{p}_r, \bar{p}_{v_0} = p^*_{v_f} + \lambda \bar {p}_v, \bar {p}_{\theta_0} = p^*_{\theta_f} + \lambda \bar {p}_\theta, \bar{p}_{\omega_0} = p^*_{\omega_f} + \lambda \bar {p}_\omega.
\label{EQ:newset}
\end{align}
where the perturbations $\lambda \bar {p}_r$, $\lambda \bar {p}_v$, $\lambda \bar {p}_\theta$ and $\lambda \bar {p}_\omega$
are uniformly selected in the intervals $[0.489,0.839]$, $[-0.317,-0.107]$, $[-0.1,0.1]$, and $[0.297,0.427]$, respectively. Calculate the value for $\bar{m}_0$ by Eq.~(\ref{EQ:Ham_para_m0}).  
Define an empty set $\mathcal D$, and insert the optimal flight-thrust vector pairs along the obtained optimal trajectories into $\mathcal D$ until the perturbation process in Eq.~(\ref{EQ:newset}) ends. Then, we can obtain a dataset $\mathcal D$ comprising the optimal state-thrust vector pairs required for NN training.



It is worth mentioning that while the values for $\bar{p}_{r_0}, \bar{p}_{v_0}, \bar{p}_{\theta_0}$, and $\bar{p}_{\omega_0}$ could be randomly chosen, the resulting value for $\bar{m}_{0}$ must satisfy $m_{d} \leq \bar {m}_{0}  \leq m_0$ during the perturbation process. Meanwhile, any flight trajectory with $\bar r(\tau) < 1, \forall~ \tau \in  [0,t_f]$ will be deleted. The propagation of the parameterized system is terminated if $\bar r(\tau) > 1.1, \forall~ \tau \in  [0,t_f]$, and $t_f$ is set as 0.9 (corresponding to an actual time of flight of $931.32$ s, which is long enough for a typical pinpoint landing).
Ultimately, a total of $26,003$ optimal trajectories are acquired, 
, and we will show in the following subsections how to use these trajectories to train NNs. 
\subsection{Using NNs to Approximate the Set-Valued Mapping $f_{\bar{S}}$}\label{setvalued}
In this subsection, we first show that $f_{\bar{S}}$ is set-valued. Then, a regularisation function will be introduced, allowing us to eliminate the one-to-many mapping.
\subsubsection{Existence of the One-to-Many Mapping}
It is common to explicitly or implicitly to assume that the input-output pairs feature a one-to-one mapping within the SL framework \citep{li2022using}. However, we will show that $f_{\bar{S}}$ is a one-to-many mapping, at least in the final phase of the powered descent, which makes it almost impossible for NNs to approximate accurately.

Substituting Eq.~(\ref{EQ:Initial_pareEquation}) into Eq.~(\ref{EQ:sf_para_system}) leads to
\begin{align}
\bar S(0) = 1 - \frac{T_{m}}{\bar {m}_0}\sqrt{\bar {p}_{v_0}^2+\bar {p}_{\omega_0}^2}.
\label{A3}
\end{align}
Combining Eq.~(\ref{EQ:Ham_para_initial_reduced}) with Eq.~(\ref{A3}), we have 
\begin{align}
\bar S(0) = \bar {p}_{v_0}.
\label{A4}
\end{align}
Recall that for the parameterized system, the pair $(\bar {p}_{v_0}, \bar {p}_{\omega_0})$ is arbitrary, and the value for $\bar {m}_0$ is determined by the pair $(\bar {p}_{v_0}, \bar {p}_{\omega_0})$ according to Eq.~(\ref{EQ:Ham_para_m0}). However, according to Eq.~(\ref{A4}), we can see that the value for $\bar S(0)$, corresponding to the flight state at touchdown, is solely affected by $\bar {p}_{v_0}$. In other words, as long as we can find different pairs $(\bar {p}_{v_0}, \bar {p}_{\omega_0})$  that will result in the same $\bar {m}_0$, the resulting value for $\bar S(0)$ will not be scalar valued.

\begin{figure}[!htp]
\begin{center}
\includegraphics[scale=0.18]{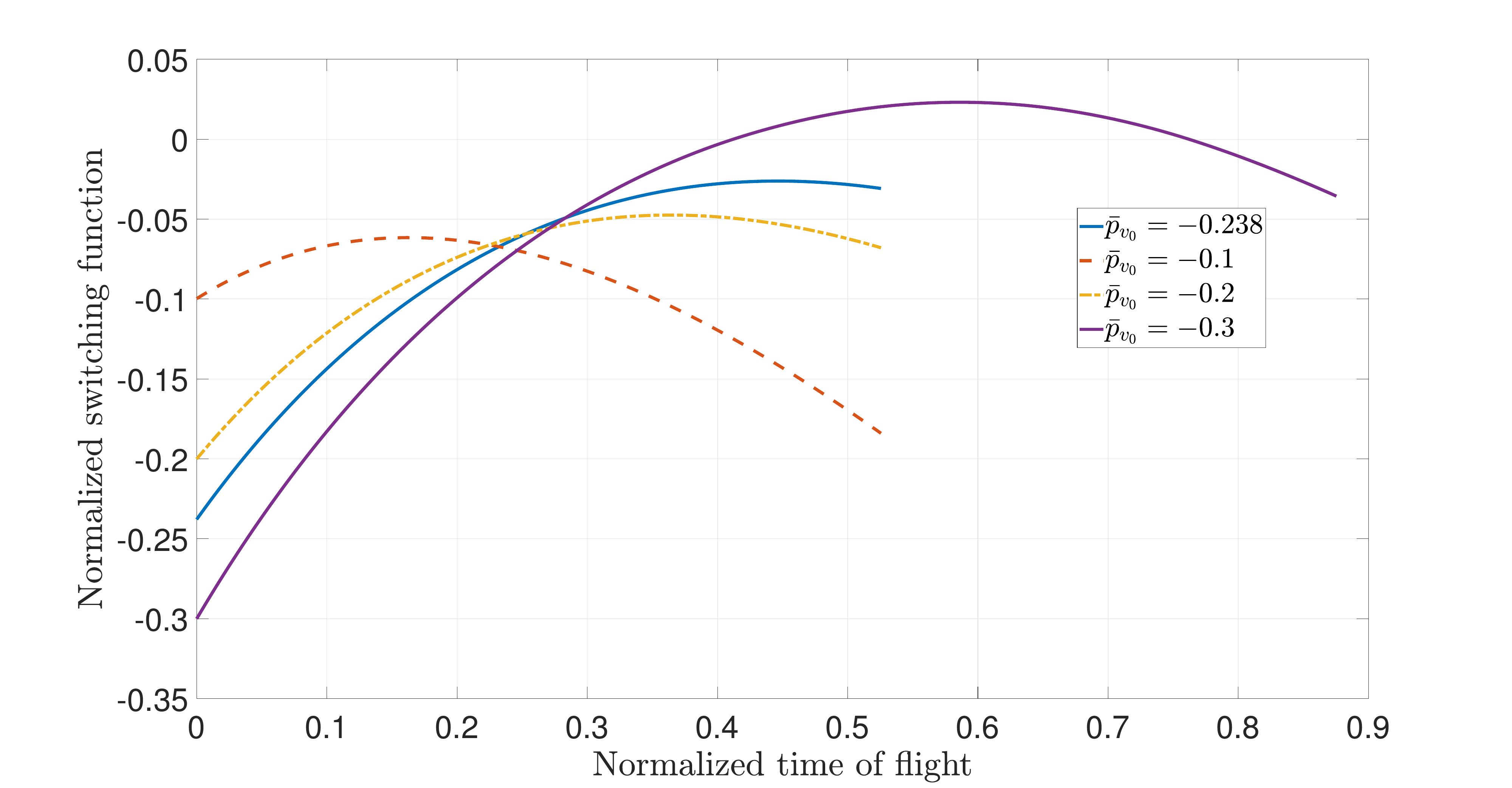}
\caption{Illustration of set-valued one-to-many mapping $f_{\bar{S}}$.}\label{Fig:onemanysf}
\end{center}
\end{figure}

For demonstration, we choose a quadruple as $(\bar {p}_{r_0},\bar {p}_{v_0}, \bar {p}_{\theta_0}, \bar {p}_{\omega_0}) = (0.753, -0.238, 0.019, 0.361)$. By solving Eq.~(\ref{EQ:Ham_para_m0}), the resulting $\bar {m}_0$ is $0.5380$. We then obtain several different values for $\bar {p}_{v_0}$, i.e., -0.1,-0.2,-0.3, and the resulting value for $\bar {m}_0$ remains unchanged. Fig.~\ref{Fig:onemanysf} shows the profiles of $\bar{S}(\tau)$ obtained by propagating the parameterized system in Eq.~(\ref{EQ:pareEquation}) with the initial condition as $\bar{\boldsymbol{x}}_{\boldsymbol{0}} = (1,0,0,0,0.5380)$.
As shown in Fig.~\ref{Fig:onemanysf}, there are different values for $\bar{S}(0)$ regarding the  same flight state at touchdown. Since the vectors $\bar {\boldsymbol x}$ and $\bar{\boldsymbol{p}}_{\boldsymbol{x}}$ are smooth everywhere, it is reasonable to conclude that at least in the final phase of the powered descent, the mapping $f_{\bar{S}}$ is also set-valued.

At touchdown, in view of Eq.~(\ref{EQ:optimal_psi_para}), it can be seen that the thrust steering angle is affected by the pair $(\bar {p}_{v_0}, \bar {p}_{\omega_0})$ used for propagating the parametrized system. Since we are able to find different pairs $(\bar {p}_{v_0}, \bar {p}_{\omega_0})$ that will result in the same flight state at touchdown. In other words, the mapping from the flight state to the thrust steering angle is also set-valued, at least in the final phase of the powered descent.  
In the case of set-valued output, the discrepancy between the prediction from the well-trained NN and the actual value in the training dataset cannot be arbitrarily reduced by simply increasing training examples or choosing different structures for the NN \citep{li2022using}.  In this case, it is impossible to use an NN to fit the mapping $f_{\bar{S}}$
, and the prediction from the well-trained NN may be incorrect or infeasible \citep{li2022using}. For this reason, we will introduce a regularisation function to eliminate the one-to-many mapping.
\subsubsection{Elimination of the One-to-Many Mapping}
\begin{figure}[!htp]
\begin{center}
\includegraphics[scale=0.18]{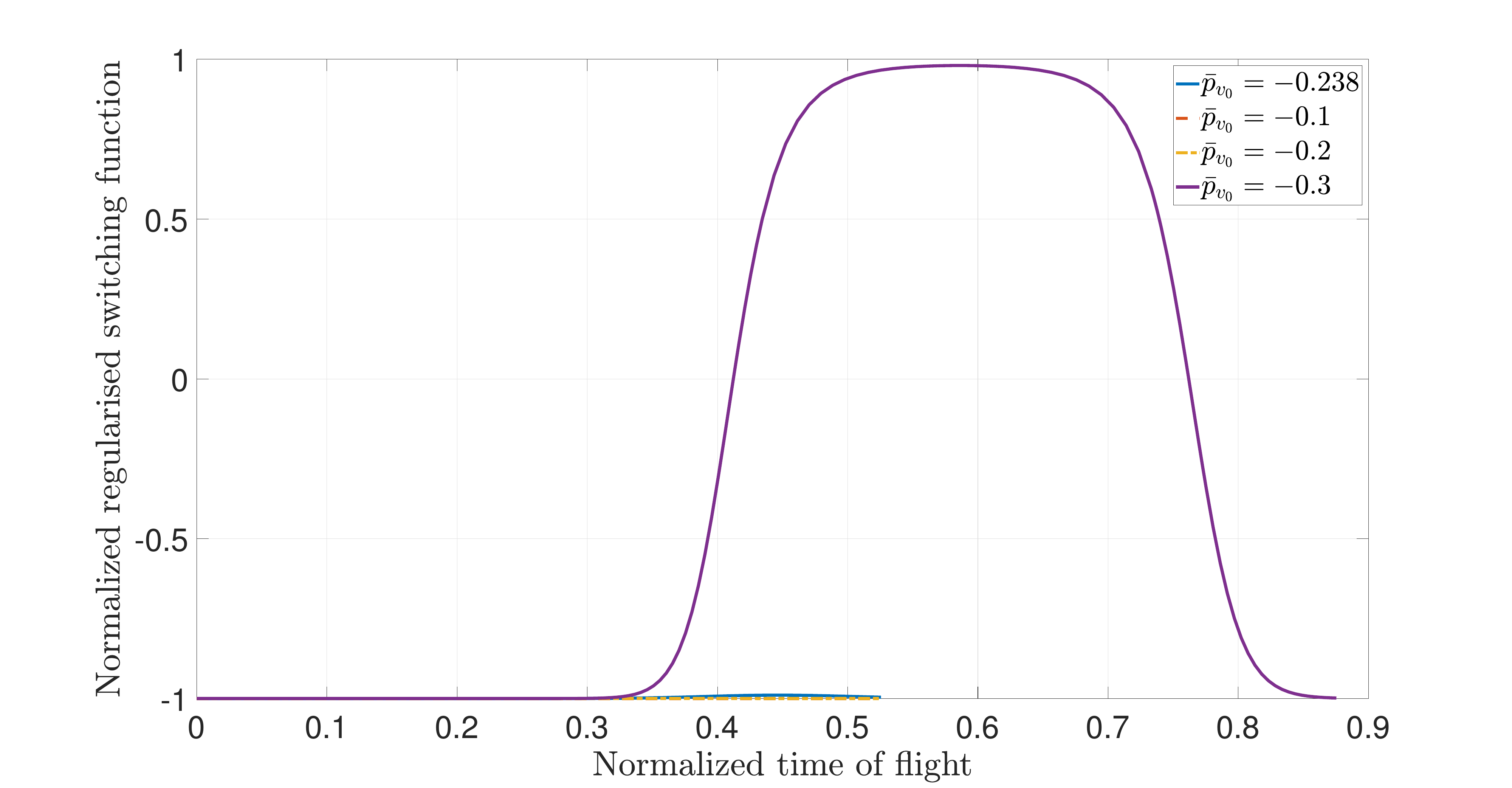}
\caption{Illustration of regularised switching functions.}\label{Fig:eliminationofonemanySF}
\end{center}
\end{figure}

For a bang-bang control profile, the switching time can be used to determine the control profile. In this regard, there are already some works on using NNs to predict the switching time \citep{wang2021real,you2021learning}. However, the number of the switching time cannot be determined \it{a prior} \rm{for the FOPDG problem}. In the final phase of the powered descent, as long as we transform the different output $\bar{S}$ into the same one, and the switching times, i.e., the zeros of the switching function, remain unchanged, the one-to-many mapping can be eliminated. For this purpose, we use a modified hyperbolic tangent function to regularise the switching function as
\begin{align}
\bar{S}_r = \tanh (\frac{\bar {S}}{\alpha}) = 1 - \frac{2}{e^{2\frac{\bar S}{\alpha}}+1},
\label{EQ:regularisation_function}
\end{align}
where $\bar {S}_r$ is the regularised switching function and $\alpha > 0$ is a small positive constant. 
There are two reasons for adopting the constant $\alpha$. Firstly, in the final phase of the powered descent, the switching function ${\bar S}$ is negative since the thrust magnitude is kept at maximum \citep{lu2023propellant}, then a small positive constant $\alpha$ will force the different switching functions to $-1$. Therefore, the mapping from the flight state to the regularised switching function will be a single value. Secondly, the constant $\alpha$ cannot be too small, or else the regularised switching function will be basically a bang-bang profile for a very small positive constant $\alpha$, which is very challenging for NNs to approximate \citep{origer2023guidance}. Thus, the constant $\alpha$ also acts like a smoothing parameter. In this paper, we set $\alpha = 0.01$, and the regularised switching function $\bar{S}_r$ related to $\bar{S}$ in Fig.~\ref{Fig:onemanysf} is displayed in Fig.~\ref{Fig:eliminationofonemanySF}. We can see that the regularised one is smooth everywhere and is not set-valued in the final phase of the powered descent, which in theory allows us to use NNs to approximate correctly \citep{origer2023guidance,li2022using}. Note that because the regularisation function does not change the zeros of the actual switching function, so there is no need for converting the regularised switching function produced by the well-trained NN back into the actual switching function.
\subsection{Scheme for Generating the FOPDG in Real Time}
With the regularisation function, we are now able to use a well-trained NN to predict the regularised switching function $\bar{S}_r$ rather than $\bar{S}$. Denote by $f_{\bar{S}_r}$ the mapping from the flight state $\boldsymbol{\bar x}$ to the regularised switching function $\bar{S}_r$. A total of three NNs are established. To elaborate, $\mathcal{D}$ is divided into three samples, i.e., $\mathcal{D}_{\tau}:=\left\{\boldsymbol{\bar x},\tau \right\}$, $\mathcal{D}_{\bar\psi}:=\left\{\boldsymbol{\bar x},\bar\psi \right\}$, and $\mathcal{D}_{\bar {S}_r}:=\left\{\boldsymbol{\bar x}, \bar {S}_r \right\}$; and the corresponding NN is denoted by $\mathcal{N}_\tau$, $\mathcal{N}_{\bar\psi}$, and $\mathcal{N}_{\bar {S}_r}$, respectively. $\mathcal{N}_{\tau}$ is designed to forecast the optimal time of flight given a flight state. 
Additionally, this network facilitates the preliminary mission design by providing time of flight evaluations without necessitating an exact solution \citep{izzo2021real}.   
$\mathcal{N}_{\bar\psi}$ predicts the thrust steering angle given a flight state. For $\mathcal{N}_{\bar {S}_r}$, it outputs the regularised switching function, which is used to determine the thrust magnitude. Once these three NNs are trained offline, they enable the real-time generation of the optimal guidance command ${\bar u}$ and $\bar{\psi}$ given a flight state.
\subsection{Training of NNs}
In this subsection, the implementation of NN training algorithm will be presented. In addition to the three aforementioned NNs, we also include the training of $\mathcal{N}_{\bar {S}}$ based on the sample $\mathcal{D}_{\bar {S}}=\left\{\boldsymbol{\bar x},\bar {S} \right\}$, in order to highlight the enhancement in approximation resulting from the regularisation function. All the networks are feedforward NNs with multiple hidden layers. Before starting the training, the dataset samples are split into $70\%$ for training, $15\%$ for validation, and $15\%$ for testing sets. Moreover, 
all input and output data are normalized by subtracting the minima and dividing by the different between the maxima and minima. For training $\mathcal{N}_\tau$, we adopt a structure with two hidden layers, each containing 15 neurons. A structure with three hidden layers, each containing 20 neurons is used to train $\mathcal{N}_{\bar\psi}$, $\mathcal{N}_{\bar {S}}$ and $\mathcal{N}_{\bar {S}_r}$.
Subsequently, the sigmoid function serves as the activation function. A linear function is utilized for the output layer. The training lies in minimizing the loss function, quantified as the Mean Squared Error (MSE) between the predicted values from the trained NNs and the actual values within the dataset samples. We employ the 'Levenberg-Marquardt' for training NNs, and the training is terminated after 1,500 epochs or when the loss function drops below $1\times 10^{-8}$. To prevent overfitting, an early stopping criteria is used. 

\begin{figure}[!htp]
\centering
\begin{subfigure}[t]{7cm}
\centering
\includegraphics[width = 8cm]{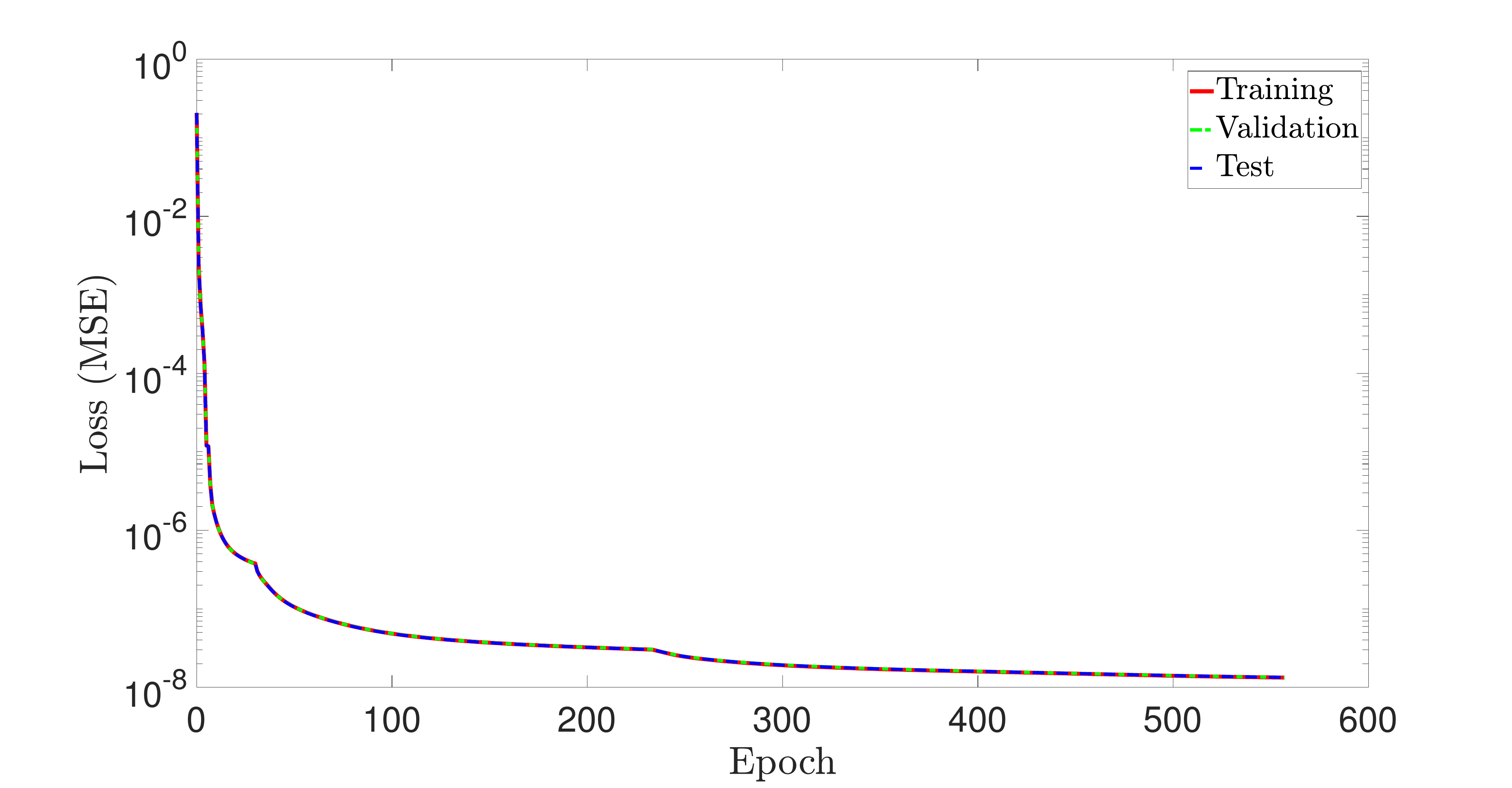}
\caption{$\mathcal{N}_{\tau}$}
\label{Fig:loss_tf}
\end{subfigure}
~~~~~
\begin{subfigure}[t]{7cm}
\centering
\includegraphics[width = 8cm]{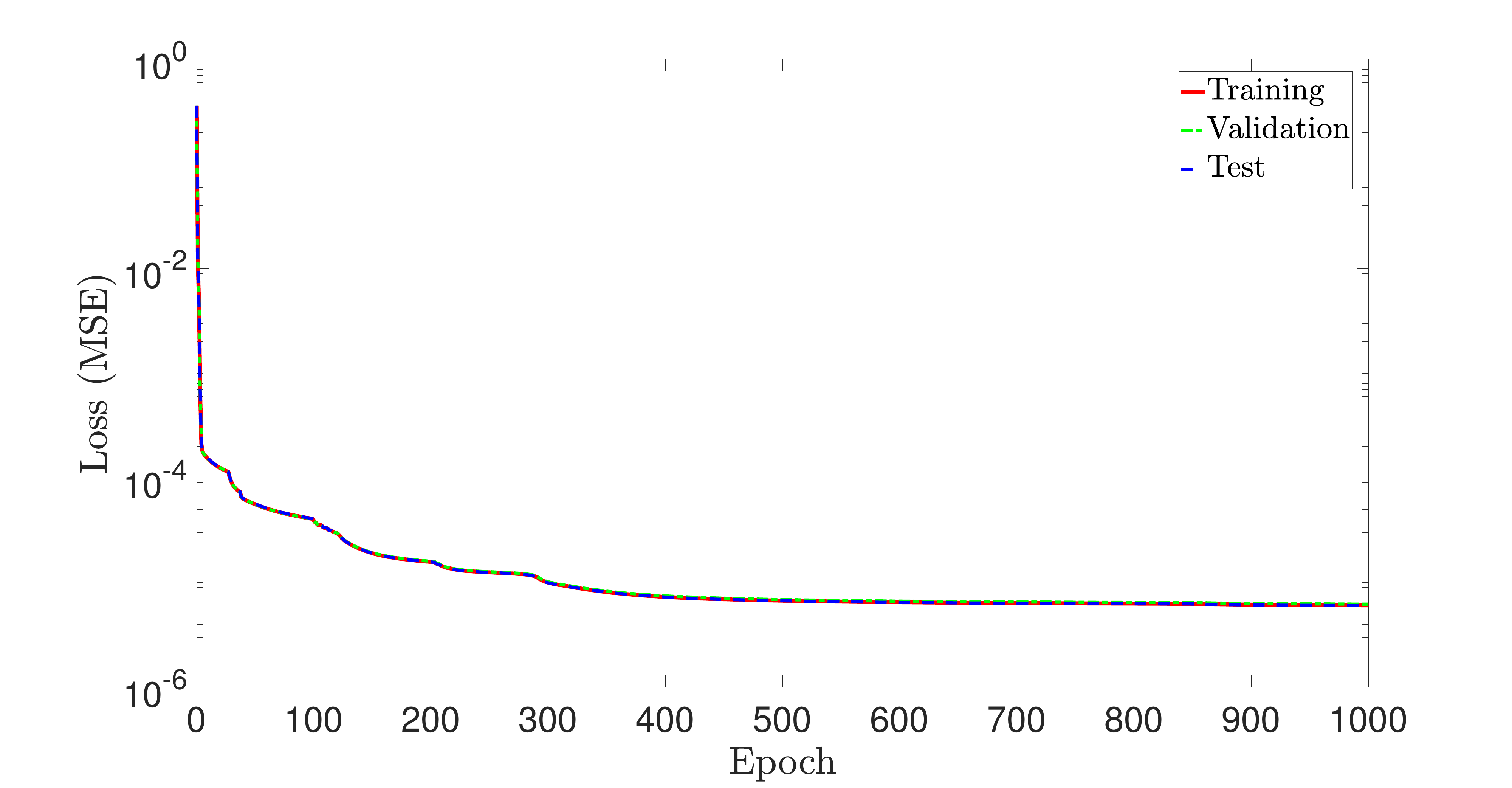}
\caption{$\mathcal{N}_{\bar\psi}$}
\label{Fig:loss_psi}
\end{subfigure}\\
\begin{subfigure}[t]{7cm}
\centering
\includegraphics[width = 8cm]{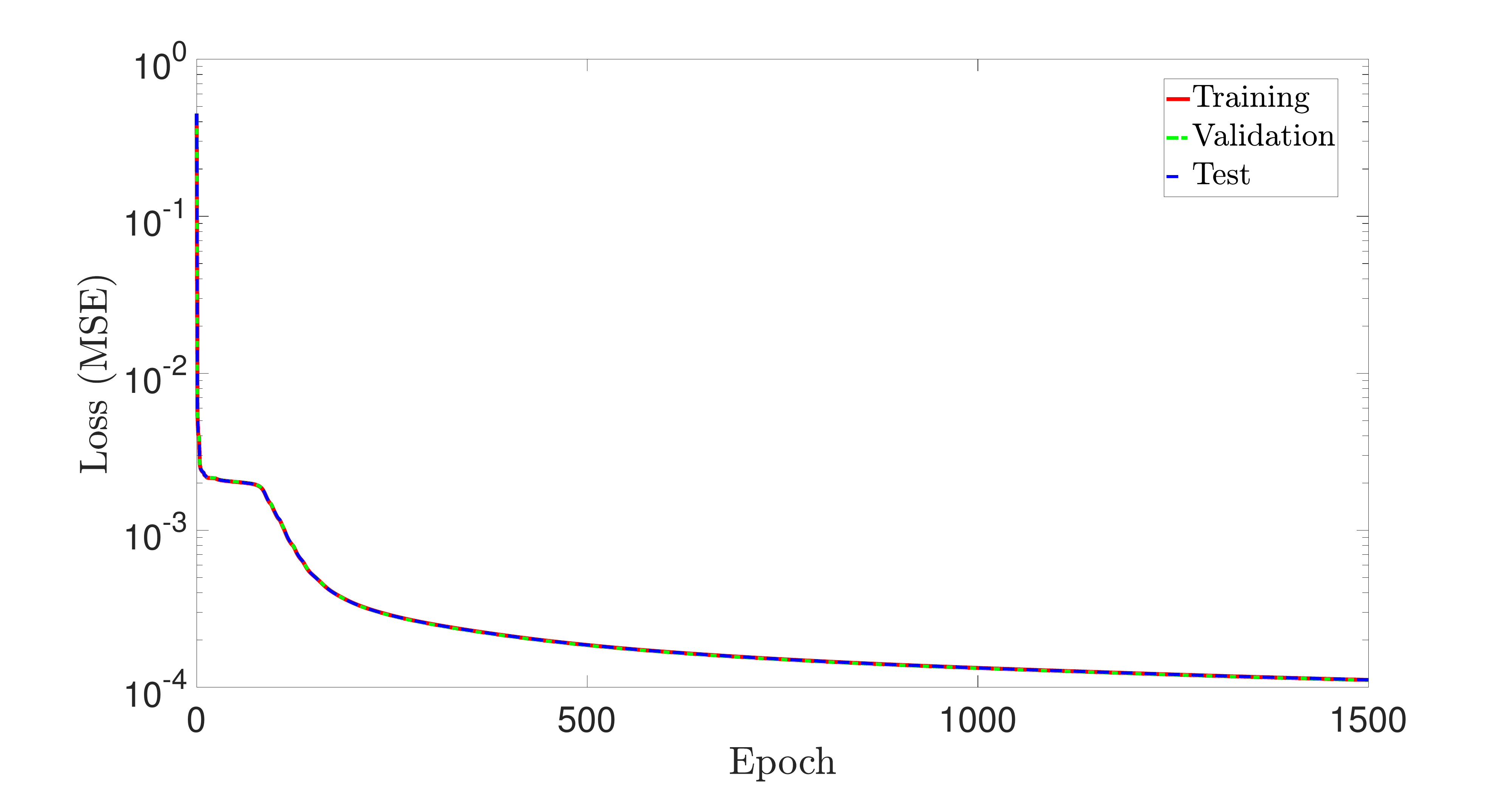}
\caption{$\mathcal{N}_{\bar {S}}$}
\label{Fig:loss_sf}
\end{subfigure}
~~~~~
\begin{subfigure}[t]{7cm}
\centering
\includegraphics[width = 8cm]{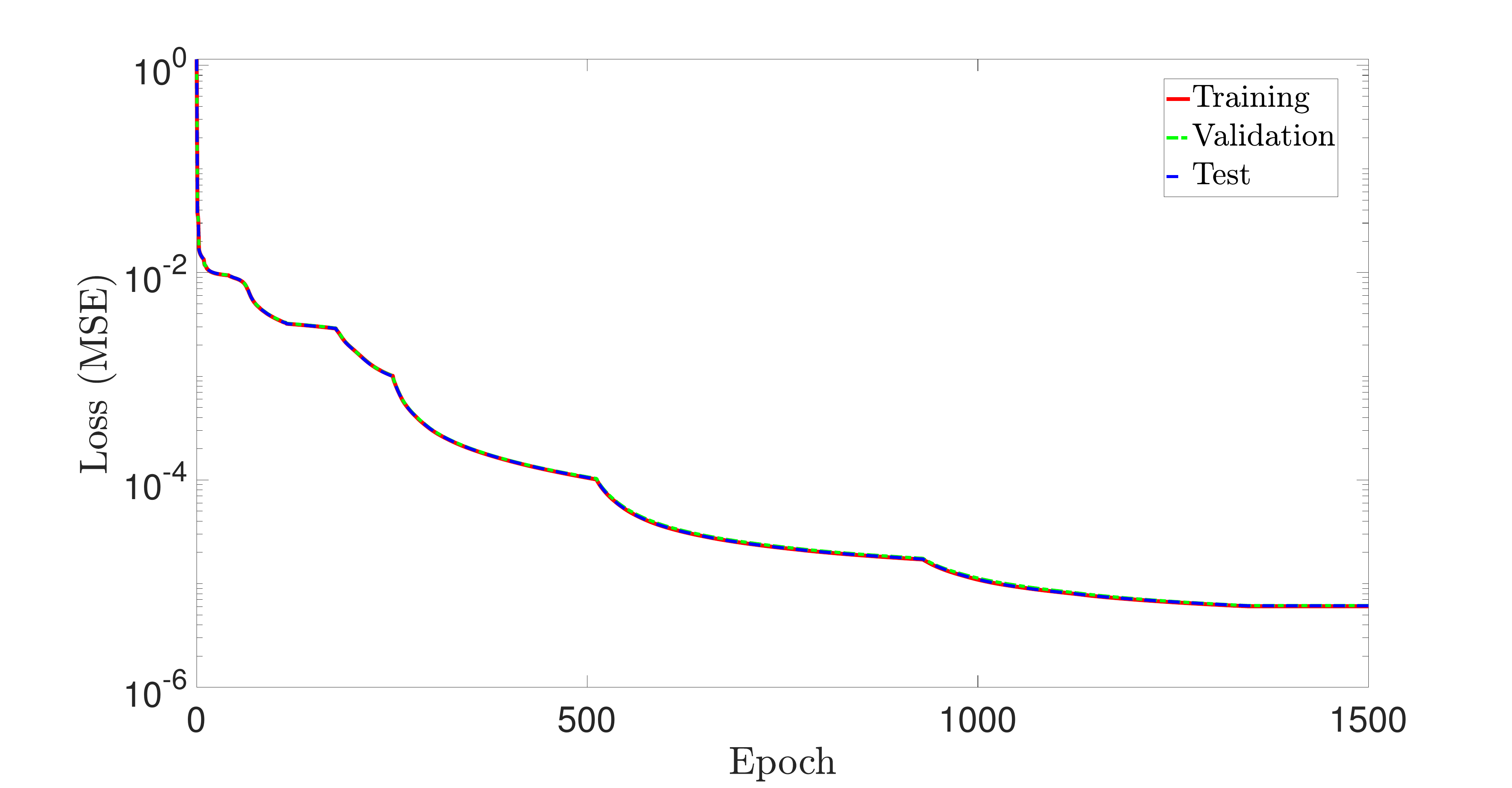}
\caption{$\mathcal{N}_{\bar {S}_r}$}
\label{Fig:loss_sf_reg}
\end{subfigure}
\caption{Training histories of four NNs.}
\label{Fig:Training}
\end{figure}

Fig.~\ref{Fig:Training} shows the training progression of the four NNs. Upon completion of the training process, the MSEs decrease to about $1.33 \times 10^{-8}$, $6.07 \times 10^{-6}$, $1.11 \times 10^{-4}$ and $6.05 \times 10^{-6}$ for $\mathcal{N}_{\tau}$, $\mathcal{N}_{\bar\psi}$, $\mathcal{N}_{\bar {S}}$, and $\mathcal{N}_{\bar {S}_r}$, respectively. Notably, the training errors of $\mathcal{N}_{\bar {S}_r}$ are greatly smaller than those of $\mathcal{N}_{\bar {S}}$, indicating that our regularisation function can generally enhance the NN's approximation accuracy. It is worth mentioning that increasing training scenarios or the size of NNs is futile for approximating a set-valued mapping \citep{li2022using}.
For further demonstration, we employ $\mathcal{N}_{\bar {S}}$ and $\mathcal{N}_{\bar {S}_r}$ to predict the switching function given the flight states from two optimal trajectories with $\bar {p}_{v_0}=-0.3$ and $\bar {p}_{v_0}=-0.2$ in Fig.~\ref{Fig:onemanysf}.

\begin{figure}[!htp]
\begin{center}
\includegraphics[scale=0.18]{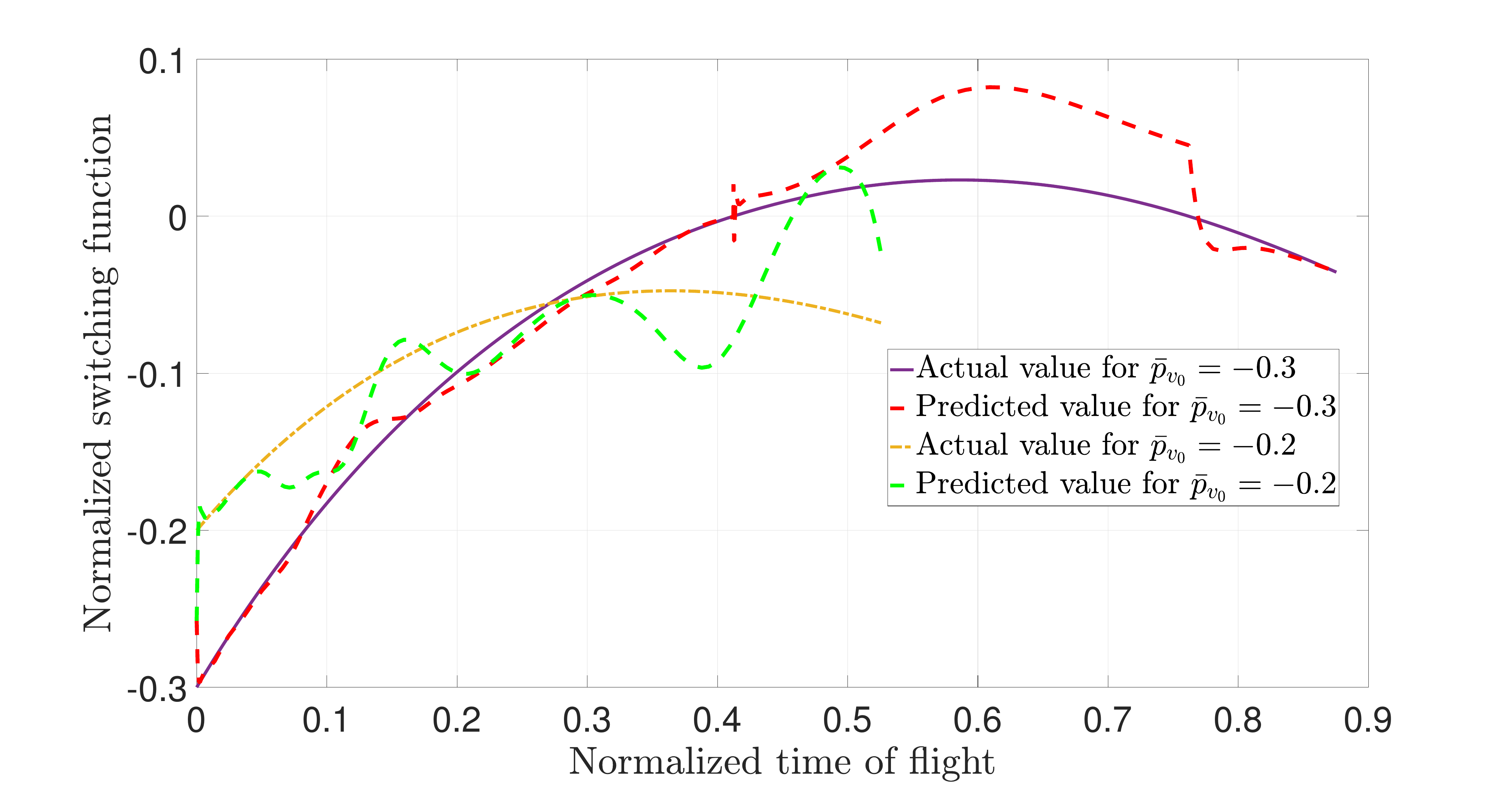}
\caption{Actual value and predicted value from $\mathcal{N}_{\bar {S}}$.}\label{Fig:conventional}
\end{center}
\end{figure}

\begin{figure}[!htp]
\begin{center}
\includegraphics[scale=0.18]{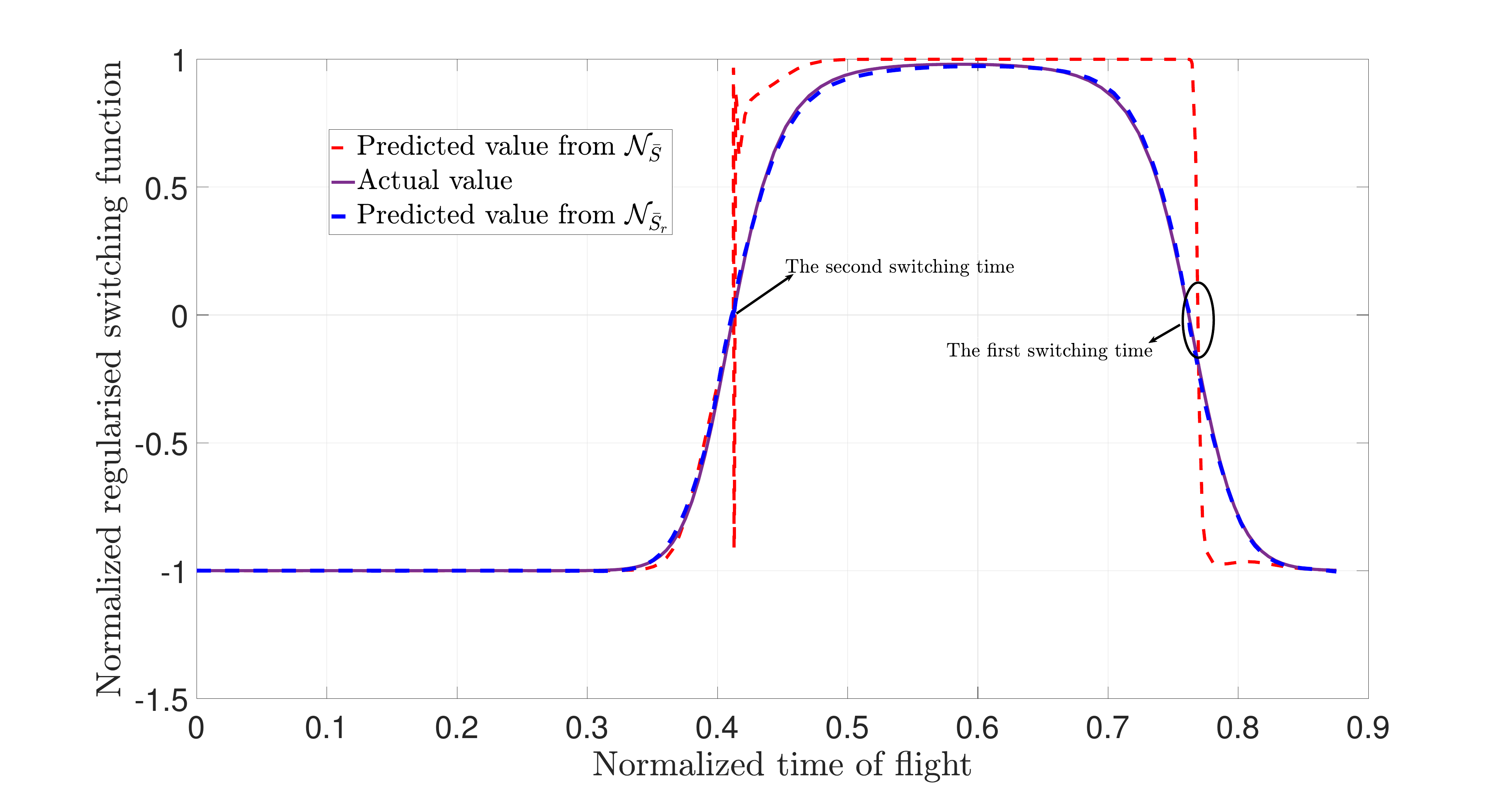}
\caption{Regularised actual and predicted values from two NNs related to the optimal trajectory with $\bar {p}_{v_0}=-0.3$.}\label{Fig:comparison_pv0}
\end{center}
\end{figure}

\begin{figure}[!htp]
\begin{center}
\includegraphics[scale=0.18]{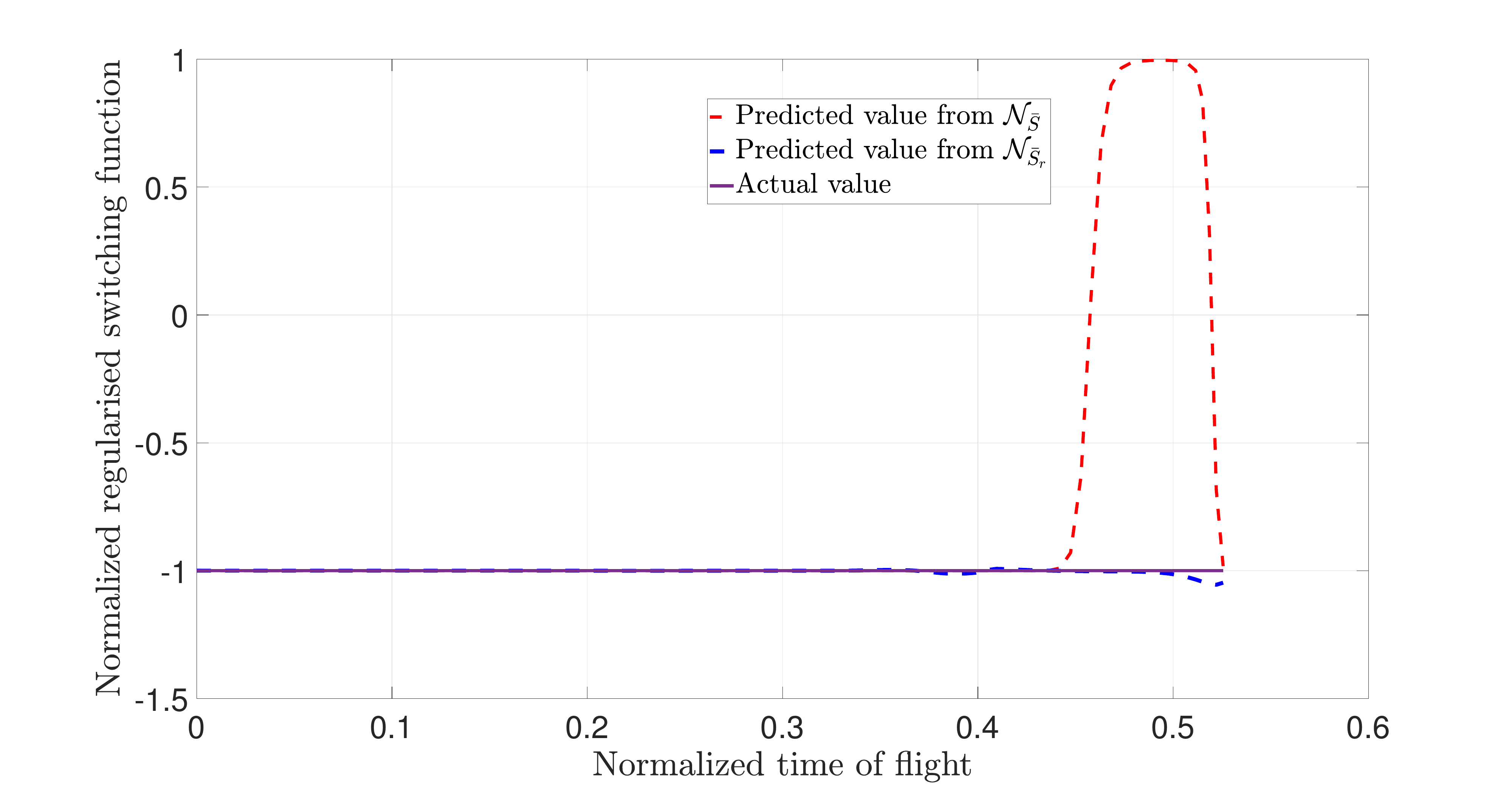}
\caption{Regularised actual and predicted values from two NNs related to the optimal trajectory with $\bar {p}_{v_0}=-0.2$.}\label{Fig:comparison_pv02}
\end{center}
\end{figure}

Fig.~\ref{Fig:conventional} displays the results of using $\mathcal{N}_{\bar S}$ to directly approximate the set-valued mapping $f_{\bar{S}}$. It is clear that the prediction from $\mathcal{N}_{\bar{S}}$ is very inaccurate for both optimal trajectories. Specifically, for the optimal trajectory with $\bar{p}_{v_0}=-0.2$, the actual value of the switching function remains negative, indicating that the optimal thrust magnitude is kept at maximum during the entire powered descent. However, the prediction from $\mathcal{N}_{\bar{S}}$ shows two zeros, which corresponds to an "on-off-on" thrust magnitude. For fair comparison, we apply the regularisation function to the actual value, as well as the predictions from $\mathcal{N}_{\bar{S}}$ and $\mathcal{N}_{\bar{S}_r}$, 
as shown in Figs.~\ref{Fig:comparison_pv0} and \ref{Fig:comparison_pv02}.
From Fig.~\ref{Fig:comparison_pv0}, we can observe that $\mathcal{N}_{\bar {S}_r}$ can accurately predict both the first and second switching time; In contrast, $\mathcal{N}_{\bar {S}}$ fails to accurately forecast the first one. As for the optimal trajectory with $\bar {p}_{v_0}=-0.2$, the predicted value from $\mathcal{N}_{\bar {S}_r}$ is quite close to the regularised real value of $-1$, as shown in Fig.~\ref{Fig:comparison_pv02}. Therefore, the predicted thrust magnitude is kept at maximum. In contrast, as expected, there is a substantial discrepancy between the prediction of  $\mathcal{N}_{\bar {S}}$ and the actual value.
\section{Numerical Simulations}\label{SE:Numerical}
To demonstrate the efficiency and effectiveness of the proposed
approach, we first compare our approach with conventional optimization methods, i.e, direct and indirect methods. Then, the landing error caused by the proposed method is investigated. All the algorithms are implemented on a desktop equipped with an Intel Core i9-10980XE CPU @3.00 GHz and 128 GB of RAM. The dynamics of the lunar lander is propagated with a maximum error of $10^{-12}$, and the flight state is updated every $0.2$ s. The simulation is terminated once the latitude of the lunar lander driven by the proposed method drops below $0.2$ m.
\subsection{Comparison with the Direct Method}
This subsection is devoted to comparing the proposed approach with an NLP solver \citep{patterson2014gpops} in terms of solution optimality and computational time. The initial condition for the lunar lander is nonnominal and listed in Table \ref{Table:cooperative1_control_effort}. 
\begin{table}[!htp]
\centering
\caption{The initial condition for comparison with the direct method}
\begin{tabular}{cccccc}
\hline
Parameter      & $r_0$~(km)     & $v_0$~(m/s)  & $\theta_0$~(deg)     & $\omega_0$~(rad/s) & $m_0$~(kg)  \\ 
\hline
Value & $1753.07$  &$-56.24$  &$4.4335$   &$5.6557 \times 10^{-4}$  &$432.44$ \\ 
\hline
\label{Table:cooperative1_control_effort}
\end{tabular}
\end{table}

\begin{figure}[!htp]
\centering
\begin{subfigure}[t]{7cm}
\centering
\includegraphics[width = 8cm]{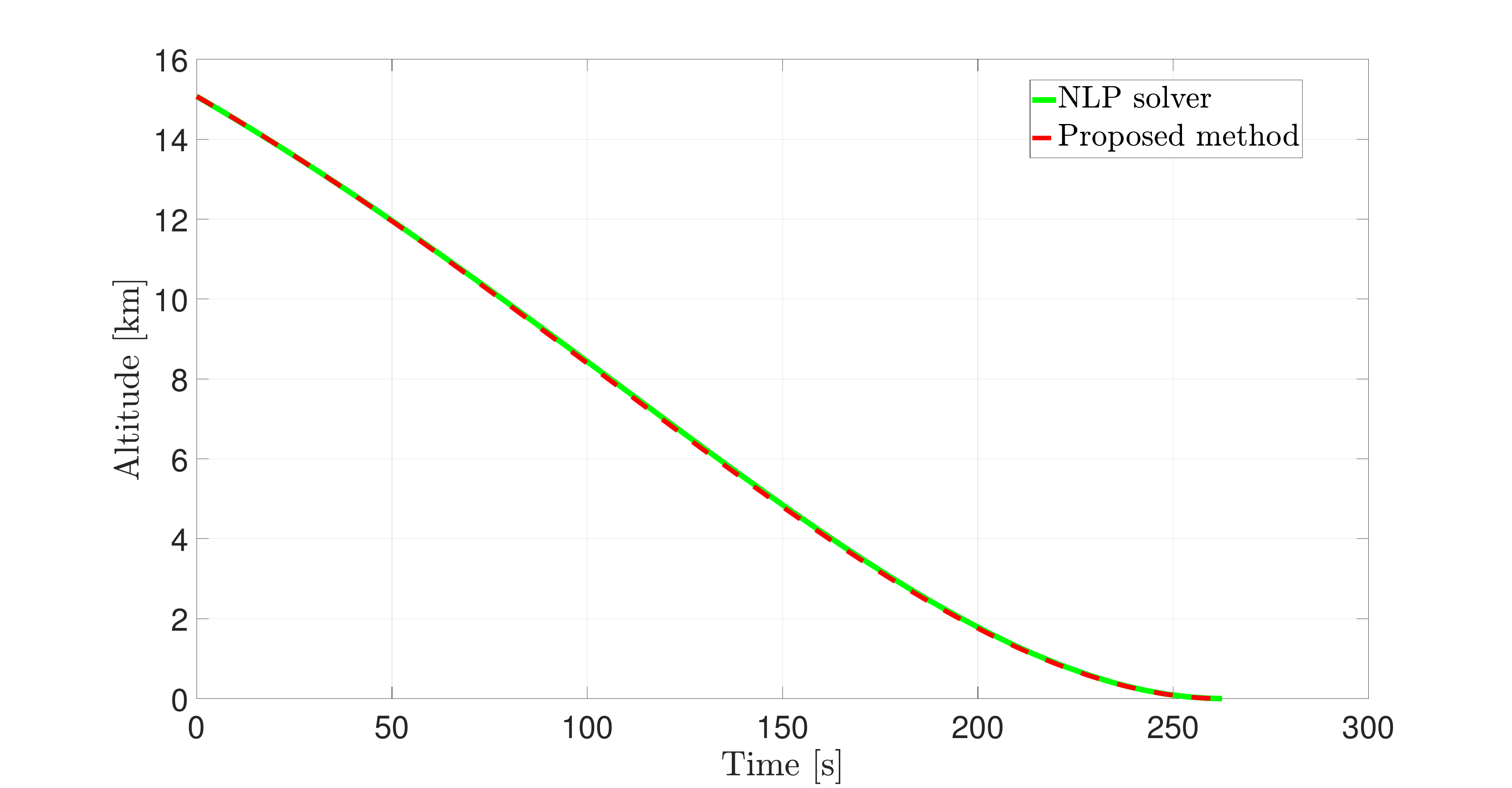}
\caption{Altitude profiles}
\label{Fig:cooperative_control_1}
\end{subfigure}
~~~~~
\begin{subfigure}[t]{7cm}
\centering
\includegraphics[width = 8cm]{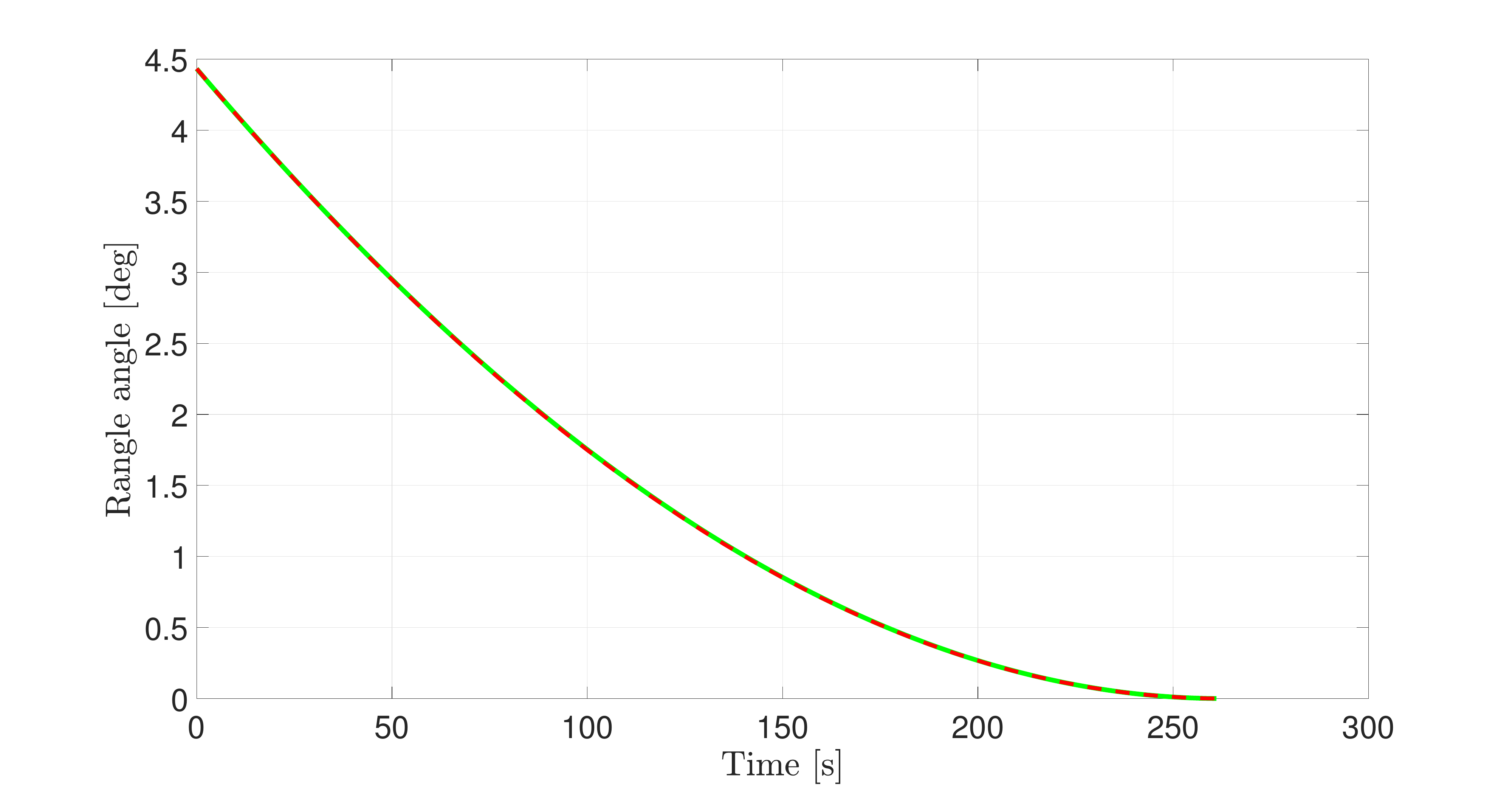}
\caption{Range angle profiles}
\label{Fig:cooperative_control_2}
\end{subfigure}\\
\begin{subfigure}[t]{7cm}
\centering
\includegraphics[width = 8cm]{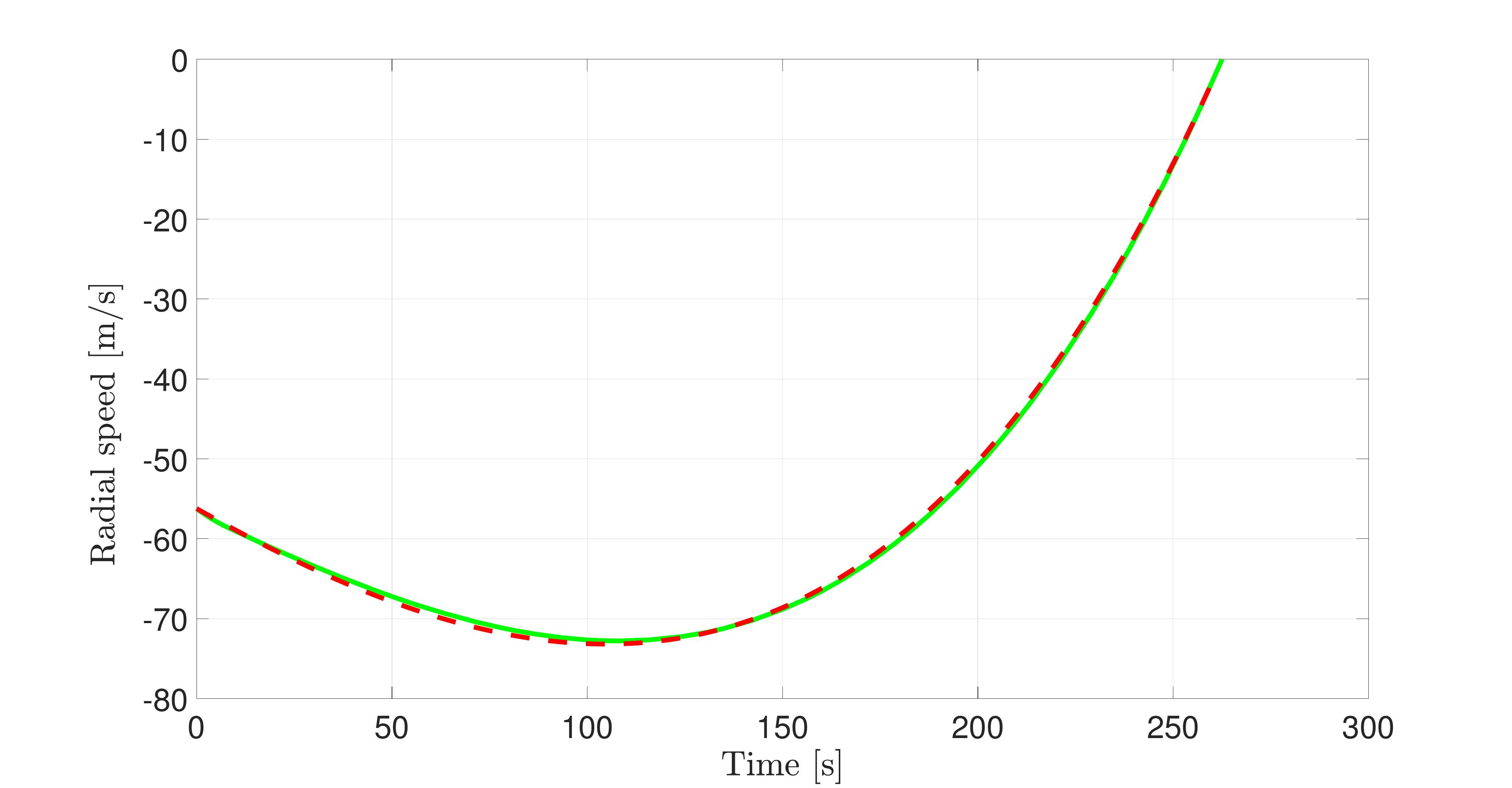}
\caption{Radial speed profiles}
\label{Fig:cooperative_control_3}
\end{subfigure}
~~~~~
\begin{subfigure}[t]{7cm}
\centering
\includegraphics[width = 8cm]{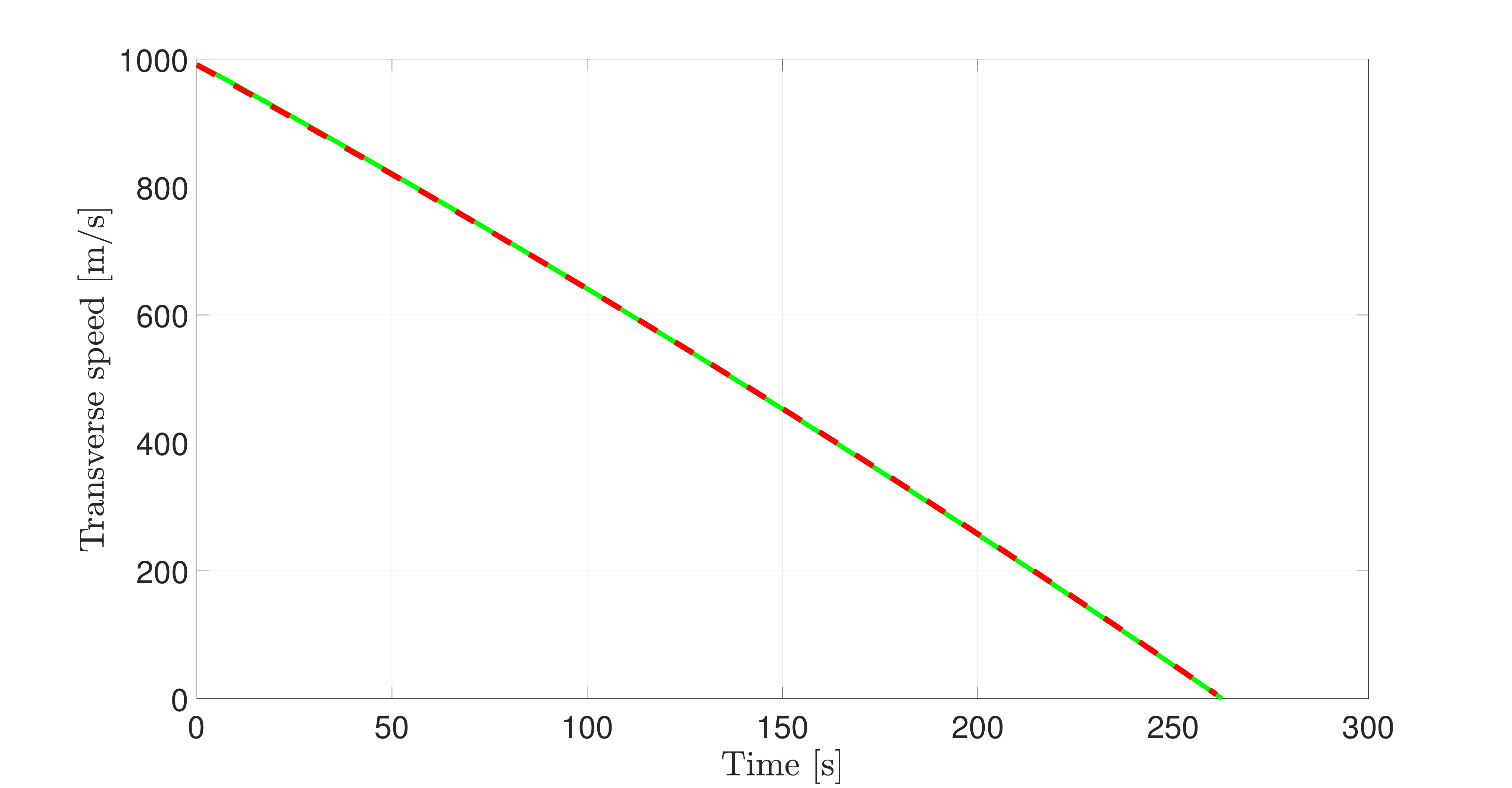}
\caption{Transverse speed profiles}
\label{Fig:cooperative_control_4}
\end{subfigure}
\caption{Profile comparisons in terms of the altitude, range angle, radial speed and transverse speed.}
\label{Fig:cooperative_control_profile}
\end{figure}

\begin{figure}[!htp]
\centering
\begin{subfigure}[t]{0.45\textwidth}
\centering
\includegraphics[scale=0.126]{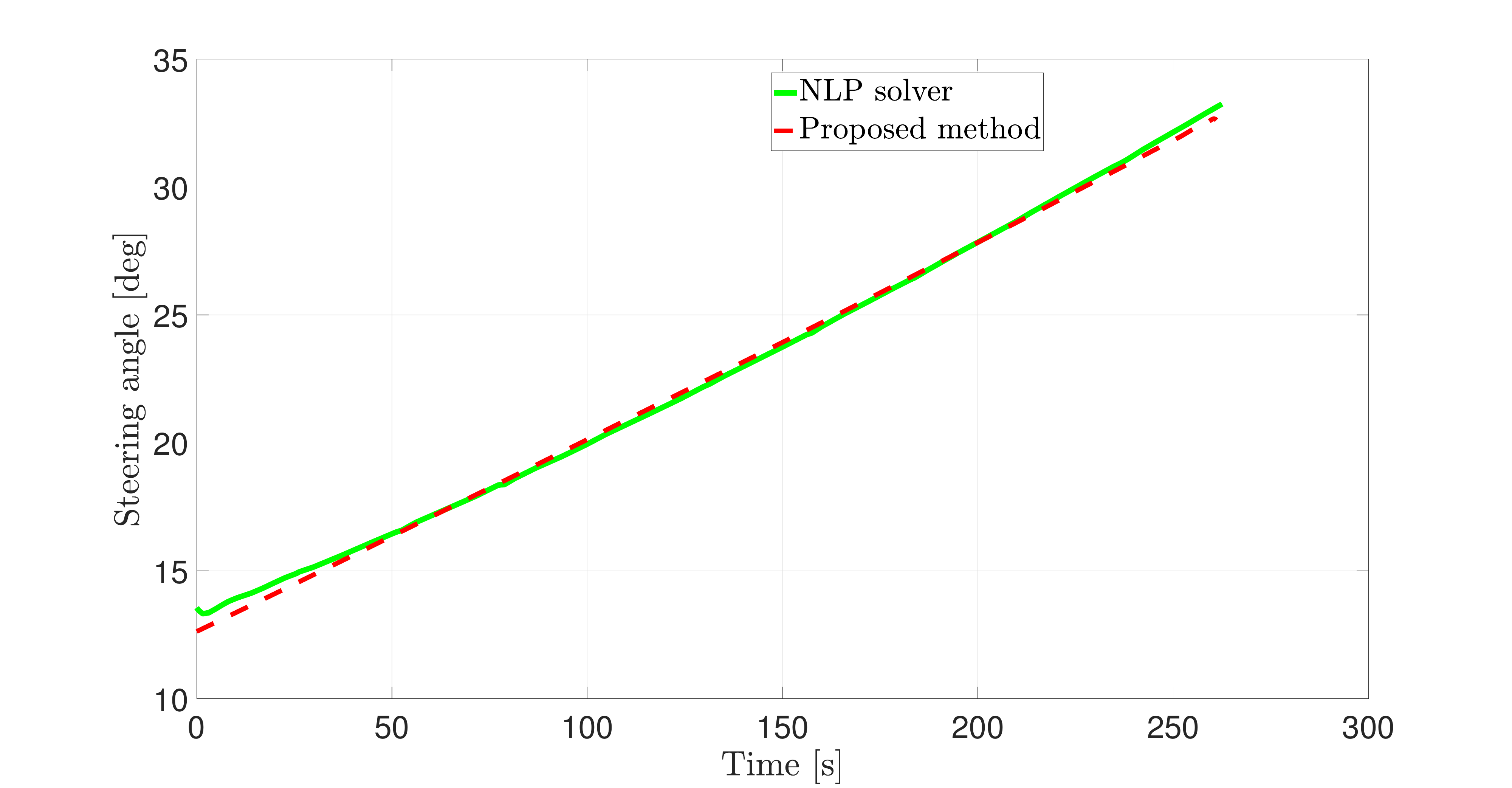}
\caption{Steering angle profiles}
\label{Fig:direction}
\end{subfigure}
~~~~~
\begin{subfigure}[t]{0.45\textwidth}
\centering
\includegraphics[scale=0.126]{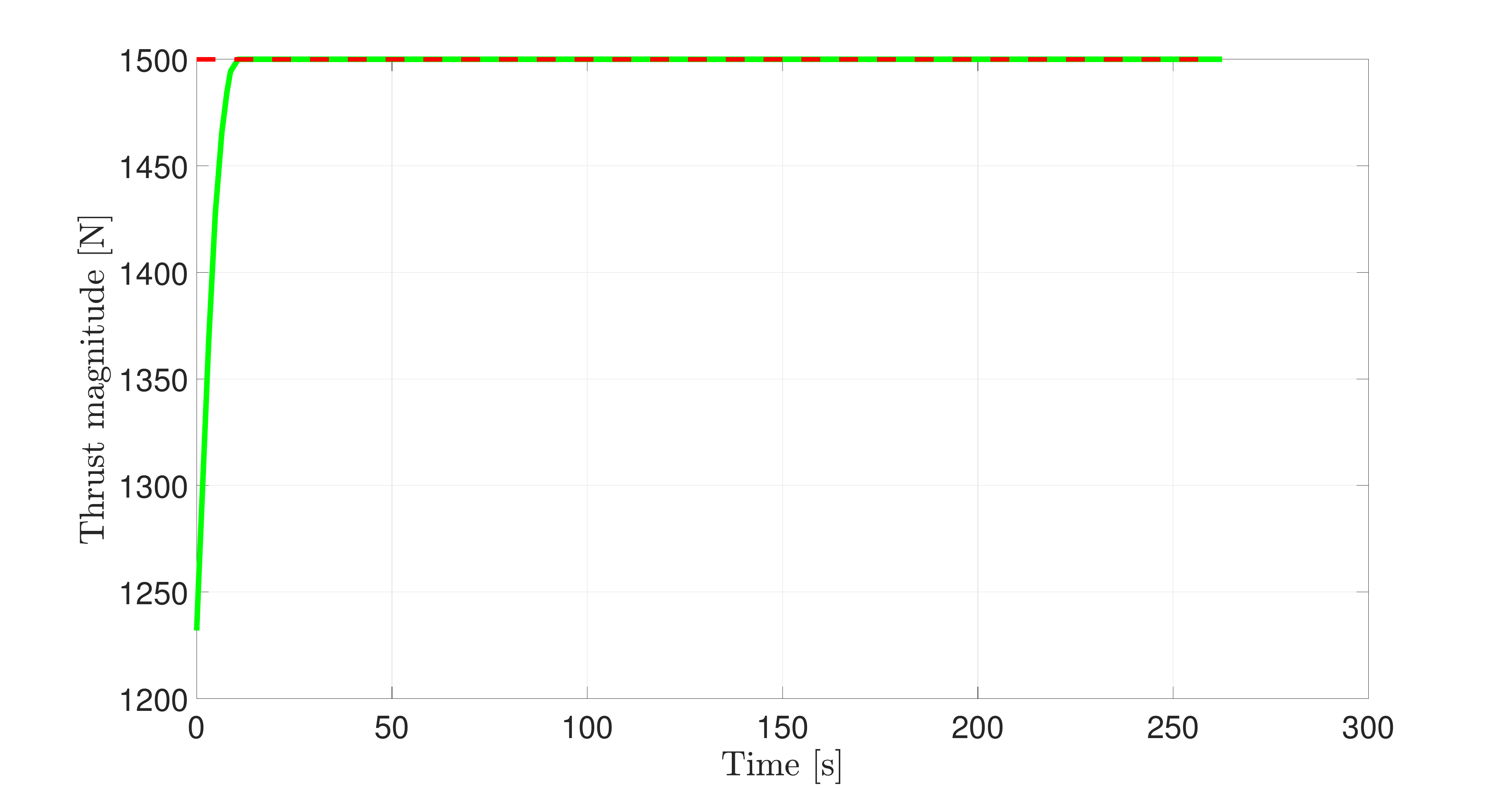}
\caption{Thrust magnitude profiles}
\label{Fig:magnitude}
\end{subfigure}\\
\caption{Thrust vector profiles obtained from the proposed method and the direct method.}
\label{Fig:NominalSolution}
\end{figure}

\begin{figure}[!htp]
\begin{center}
\includegraphics[scale=0.18]{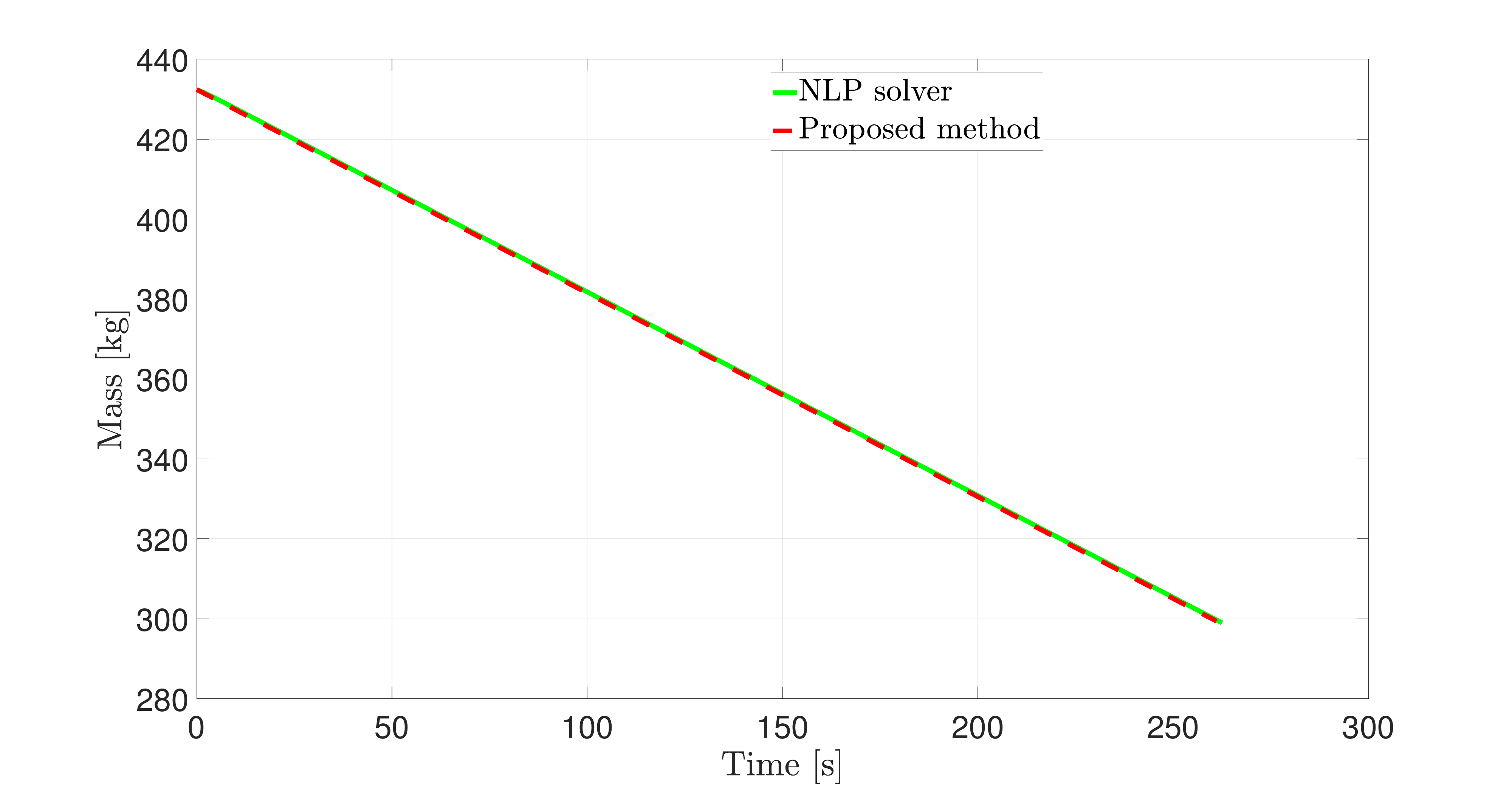}
\caption{Mass profiles obtained from the proposed method and the direct method.}\label{Fig:Mass_case1}
\end{center}
\end{figure}

Fig.~\ref{Fig:cooperative_control_profile} shows some state profiles obtained from the proposed method and the NLP solver, from which we can see that both methods guide the lunar lander to the desired landing site successfully. Specifically,  the transverse speed profiles generated by the two methods 
 decrease almost linearly to near zero, as shown in Fig.~\ref{Fig:cooperative_control_4}. 
The thrust vector profiles obtained from these two methods are displayed in Fig.~\ref{Fig:NominalSolution}, 
from which we can see that both the thrust steering angle and thrust magnitude are different. Notably, the direct method fails to generate a bang-bang thrust magnitude, while the thrust magnitude obtained from the proposed method is kept at maximum during the entire powered descent. It is worth mentioning that the thrust magnitude generated by the direct method may not be flyable by a real engine
without the capability of thrust magnitude regulation.
In addition, the mass profiles obtained via the two methods are displayed in Fig.~\ref{Fig:Mass_case1}.  
As a result, it takes $261.0667$ s for the lunar lander to arrive at the landing site with a fuel consumption of $133.06$ kg by adopting the proposed method. In contrast, the final time obtained by the NLP solver is $262.5044$ s, and the corresponding fuel consumption is $133.50$ kg. 

To assess the performance of fulfilling the requirement for nearly zero touchdown speed, we denote by $V_f$ the terminal speed error at touchdown caused by the proposed method, and it can be determined as
\begin{align*}
V_f = \left\|(v^{\mathcal{N}}_f, \omega^{\mathcal{N}}_f * r^{\mathcal{N}}_f )\right\|_2, 
\end{align*}
where $v^{\mathcal{N}}_f$, $\omega^{\mathcal{N}}_f$, and $r^{\mathcal{N}}_f$ specifies the terminal radial speed, angular velocity, and radial distance of the lunar lander, respectively.  As a result, the terminal speed error of the lunar lander caused by the proposed method is $0.4268$ m/s, which is small enough to avoid crashing the lunar lander.

Now we compare the computational time. Recall that our method requires two well-trained NNs  to generate the thrust vector (${\mathcal{N}}_{\tau}$ is only used to predict the time of flight, and it is not used for generating the guidance command). A total of 10,000 trials of the proposed method across various flight states are run in a C-based computational environment, and the mean execution time for generating the thrust vector is 0.0026 ms. This translates to approximately 0.0780 ms on a typical flight processor operating at 100 MHz \citep{gankidi2017fpga}. On the other hand, the NLP solver takes around $3.3752$ s to find the thrust vector, in which the thrust magnitude may not be exactly bang-bang. 
\subsection{Comparison with the Indirect Method}
Since the direct method may fail to generate a bang-bang thrust magnitude, we consider to compare the proposed method with the indirect method, which resolves the shooting function in Eq.~(\ref{EQ:TPBVP_law}).  The nonnominal initial condition of the lunar lander is presented in Table \ref{Table:cooperative1_control_effort_shooting}.
\begin{table}[!htp]
\centering
\caption{The initial condition for comparison with the indirect method}
\begin{tabular}{cccccc}
\hline
Parameter      & $r_0$~(km)     & $v_0$~(m/s)  & $\theta_0$~(deg)     & $\omega_0$~(rad/s) & $m_0$~(kg)  \\ 
\hline
Value & $1762.05$  &$21.35$  &$24.02$   &$1.1274 \times 10^{-3}$  &$600$ \\ 
\hline
\label{Table:cooperative1_control_effort_shooting}
\end{tabular}
\end{table}

\begin{figure}[!htp]
\begin{center}
\includegraphics[scale=0.18]{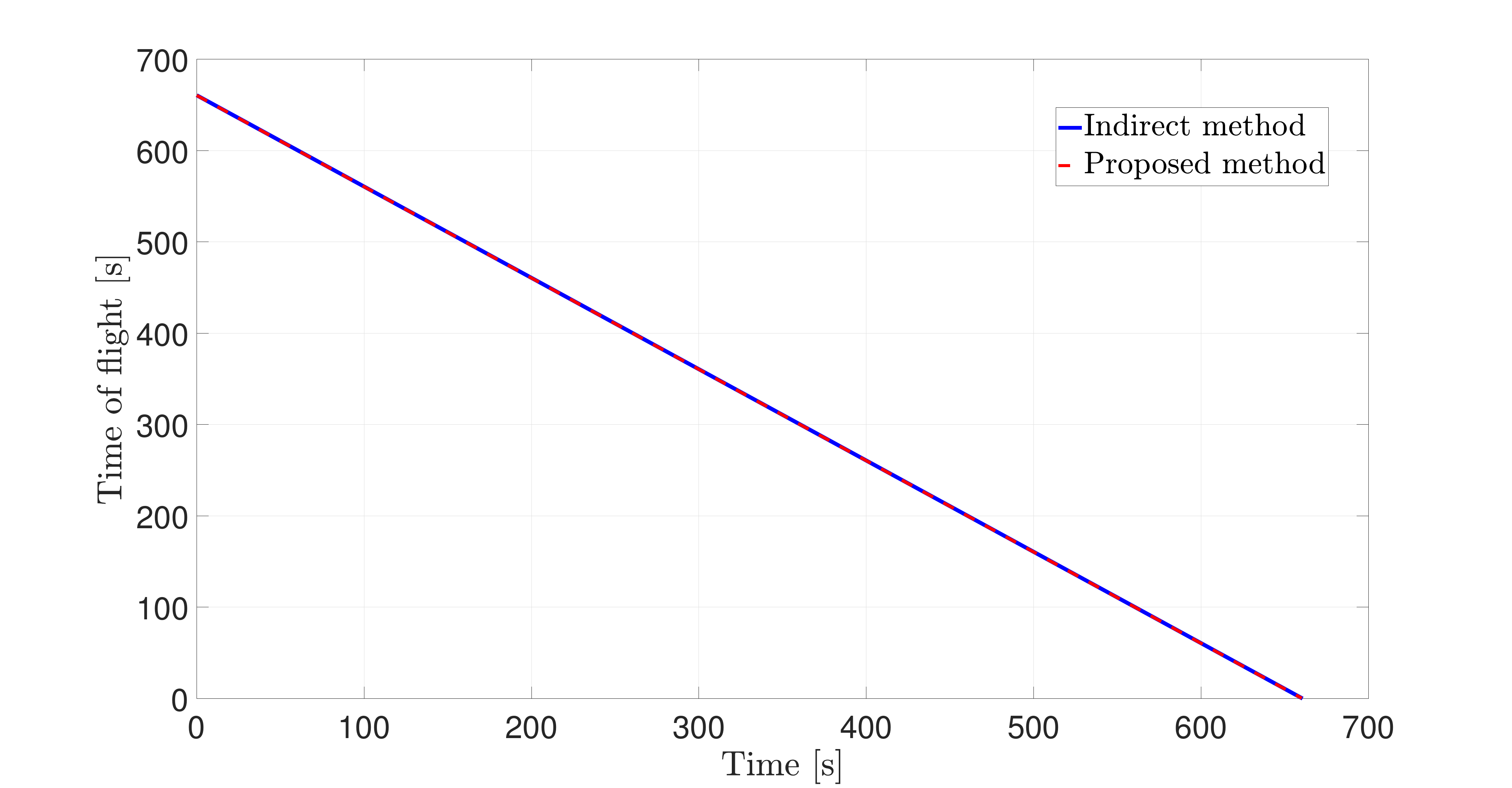}
\caption{Time of flight profiles obtained from the proposed method and indirect method.}\label{Fig:timeofflight_comparsion}
\end{center}
\end{figure}

\begin{figure}[!htp]
\centering
\begin{subfigure}[t]{7cm}
\centering
\includegraphics[width = 8cm]{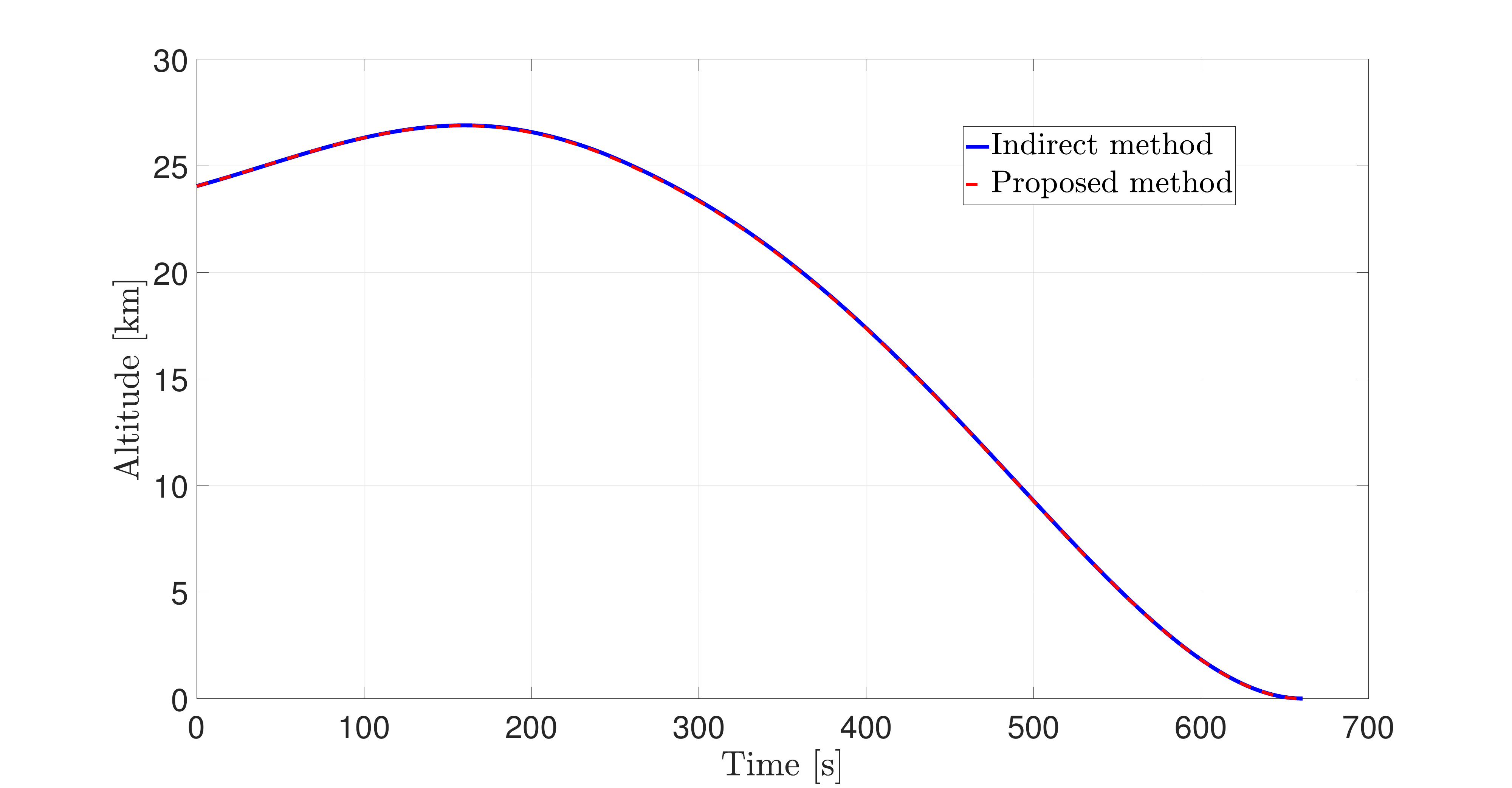}
\caption{Altitude profiles}
\label{Fig:cooperative_control_1_shooting}
\end{subfigure}
~~~~~
\begin{subfigure}[t]{7cm}
\centering
\includegraphics[width = 8cm]{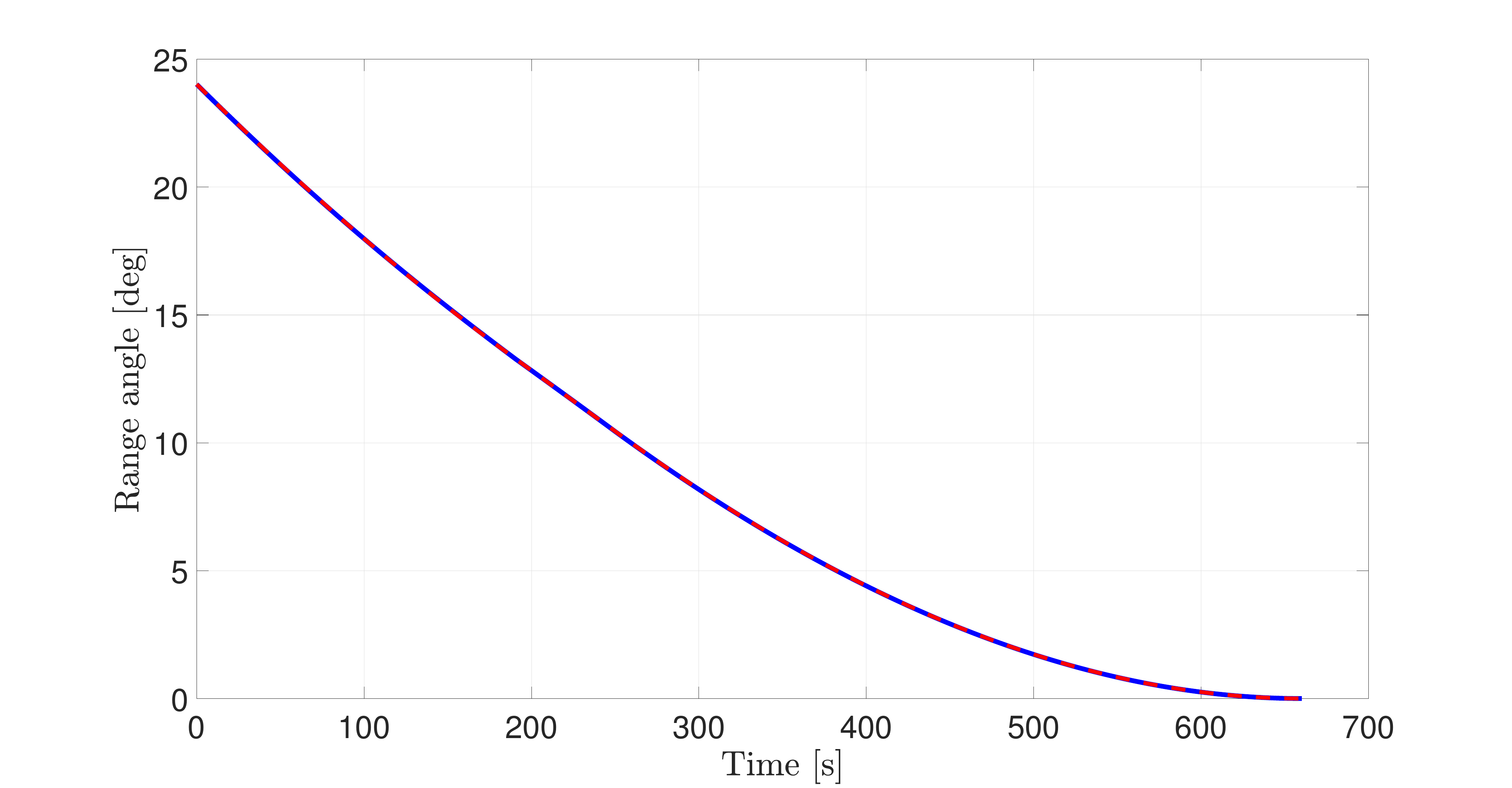}
\caption{Range angle profiles}
\label{Fig:cooperative_control_2_shooting}
\end{subfigure}\\
\begin{subfigure}[t]{7cm}
\centering
\includegraphics[width = 8cm]{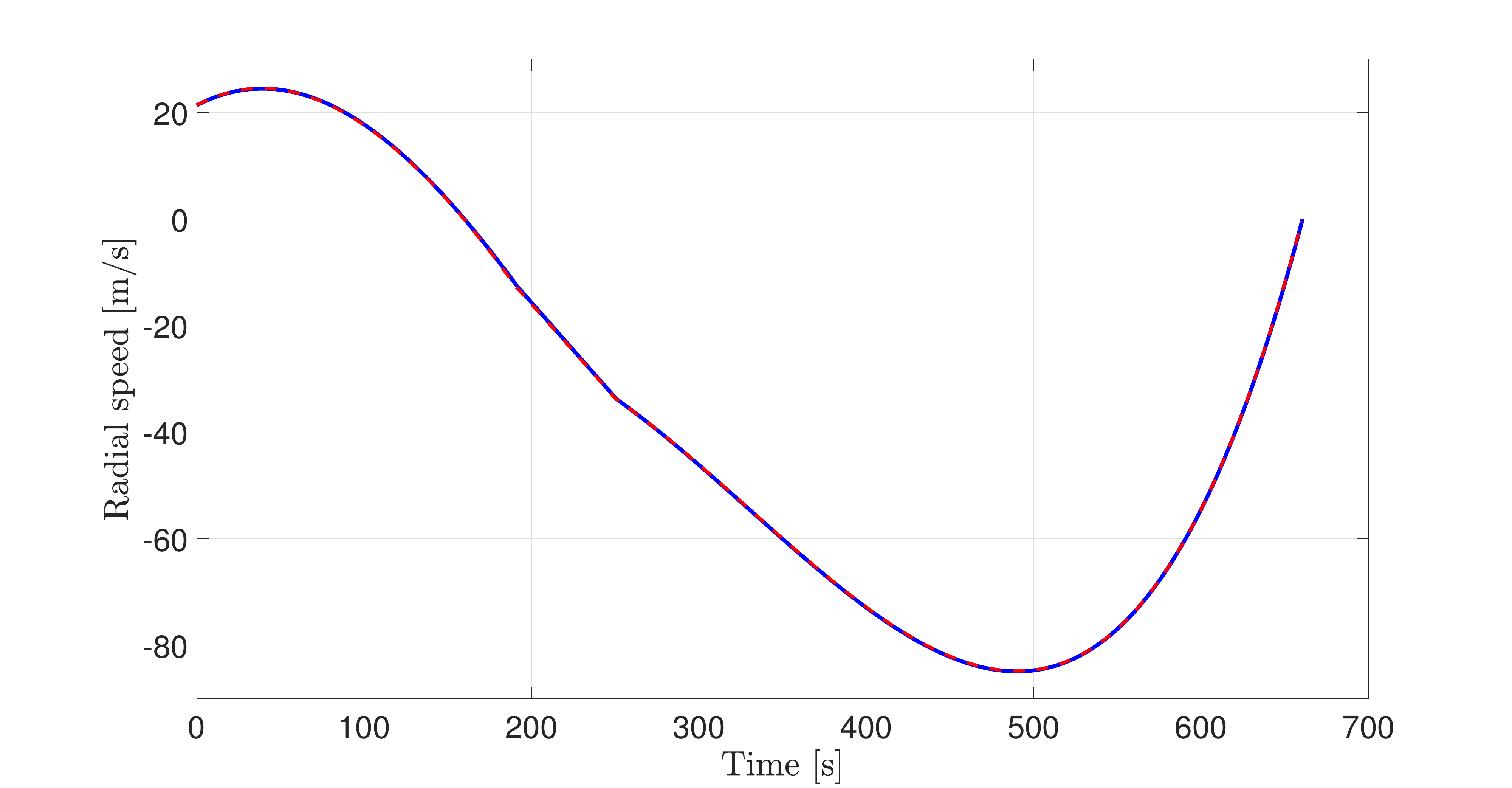}
\caption{Radial speed profiles}
\label{Fig:cooperative_control_3_shooting}
\end{subfigure}
~~~~~
\begin{subfigure}[t]{7cm}
\centering
\includegraphics[width = 8cm]{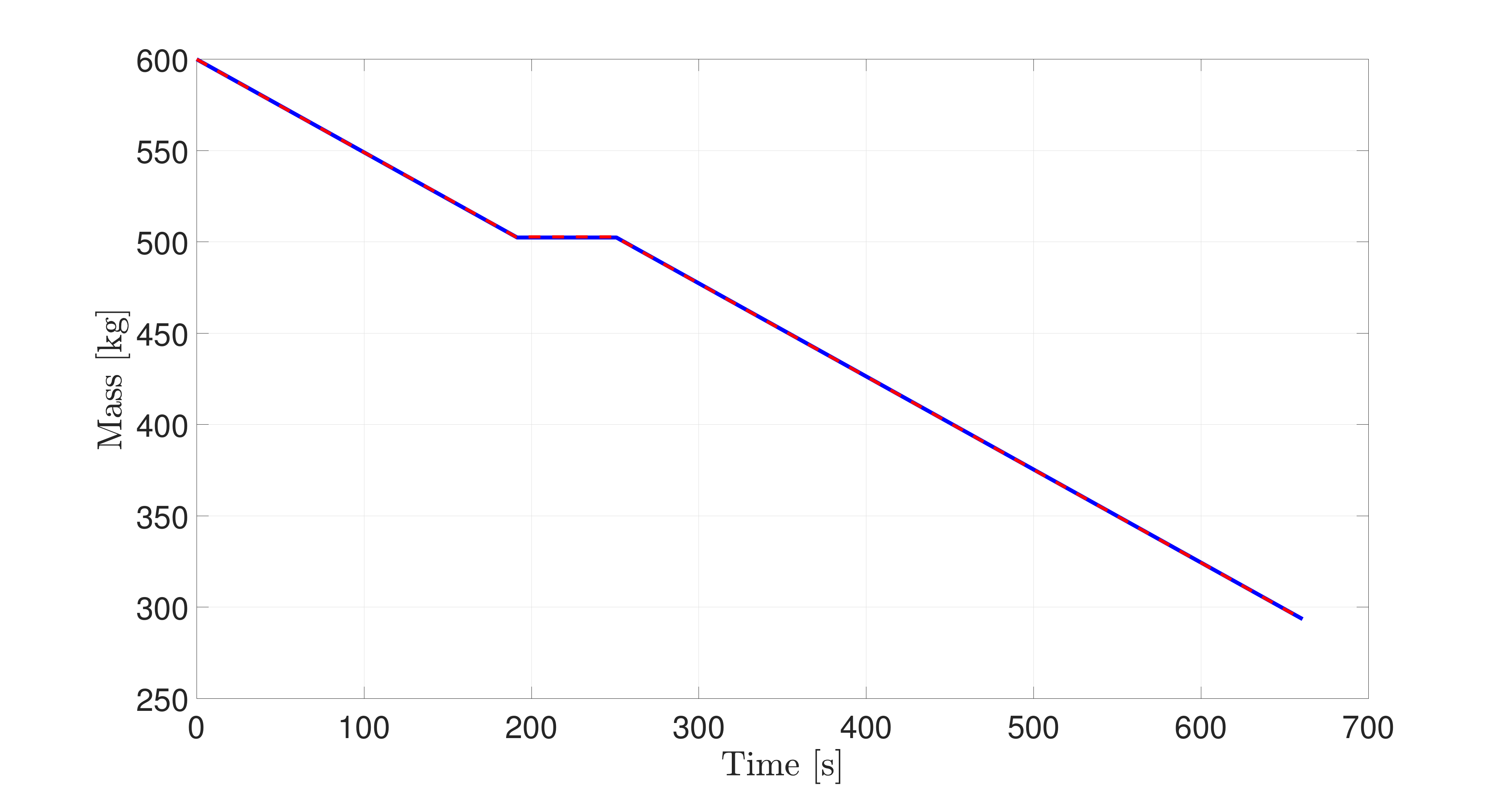}
\caption{Mass profiles}
\label{Fig:cooperative_control_4_shooting}
\end{subfigure}
\caption{Profile comparisons in terms of the altitude, range angle, radial speed and mass.}
\label{Fig:comparison_shooting}
\end{figure}

\begin{figure}[!htp]
\centering
\begin{subfigure}[t]{0.45\textwidth}
\centering
\includegraphics[scale=0.126]{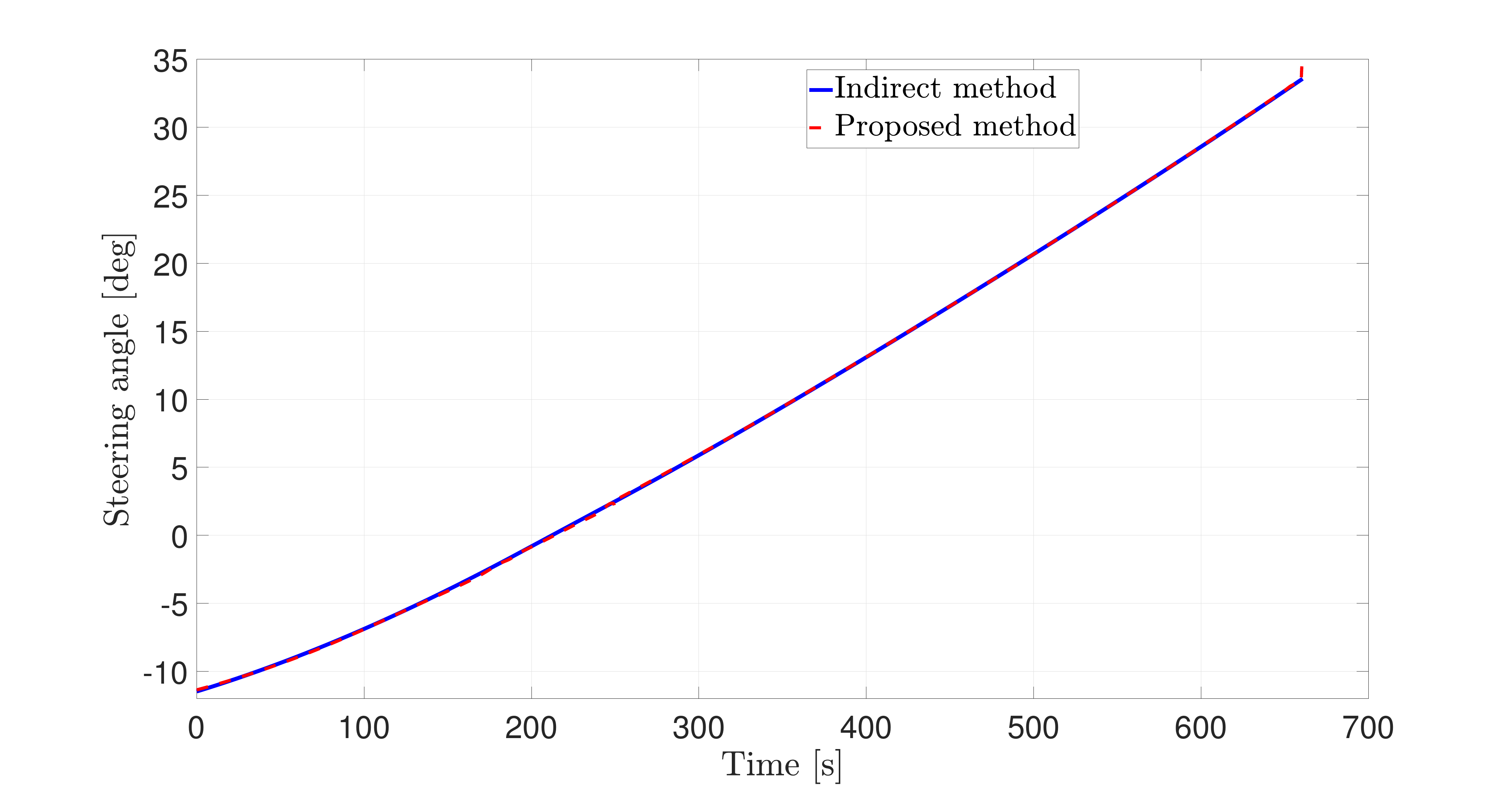}
\caption{Steering angle profiles}
\label{Fig:direction_shooting}
\end{subfigure}
~~~~~
\begin{subfigure}[t]{0.45\textwidth}
\centering
\includegraphics[scale=0.126]{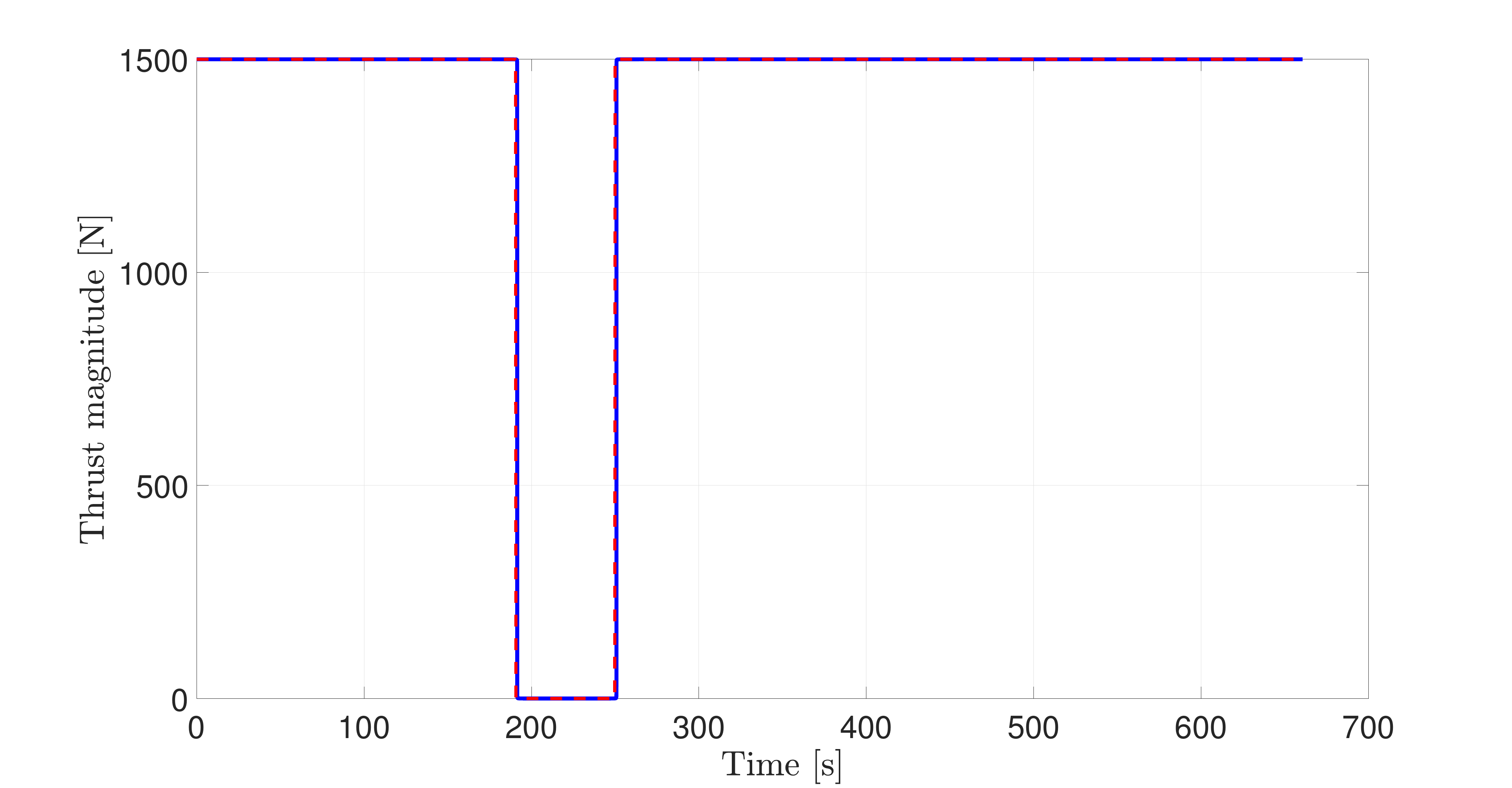}
\caption{Thrust magnitude profiles}
\label{Fig:magnitude_shooting}
\end{subfigure}\\
\caption{Thrust vector profiles obtained from the proposed method and indirect method.}
\label{Fig:ShootingSolution}
\end{figure}

The prediction of the time of flight for the given initial condition is $660.48$ s, while the final time obtained from the indirect method is $660.62$ s. Fig.~\ref{Fig:timeofflight_comparsion} presents the time of flight profiles obtained from the proposed method and indirect method during the powered descent, from which we can see that, ${\mathcal{N}}_{\tau}$ is able to precisely predict the time of flight during the entire powered descent. 
Fig.~\ref{Fig:comparison_shooting} compares the solution obtained from the two methods in terms of the altitude, range angle, radial speed, and mass. It can be seen that the solutions are almost identical to each other, which demonstrates high accuracy of the proposed method. 
In addition, the thrust vector profiles via the two different methods are shown in Fig.~\ref{Fig:ShootingSolution}. The thrust steering angle profile obtained from the proposed method almost coincide with the solution from the indirect method, except at the end of the powered descent, which further verifies the fact that the mapping from the flight state to the thrust steering angle is also set-valued in such phase. Moreover, the thrust magnitude profile obtained from the proposed method takes a bang-bang control form with two times of switch, which exactly matches the solution via the indirect method, as shown in Fig.~\ref{Fig:magnitude_shooting}. The lunar lander guided by the proposed method consumes a total of $306.80$ kg fuel, and lands at the desired landing site with a terminal speed error $V_f$ of $1.4359$ m/s. In contrast, the lunar lander guided by the indirect method consumes a total of $306.49$ kg fuel during the powered descent. In terms of the computational time, a homotopy on the smoothing constant is usually required to facilitate the convergence of the indirect method. As a result, it normally takes $6.3020$ s for the indirect method to find the optimal solution. Conversely, our method is able to generate the optimal thrust vector in 0.0026 ms with a small penalty on the fuel consumption.
\subsection{Landing Error Analysis}
In the last two subsections, we have shown that the proposed method is able to generate the FOPDG in real time. However, the ability of the proposed method to accomplish the pinpoint landing merits further studying since small errors could be propagated through the flight trajectory, which may result in suboptimal or failed landings \citep{sanchez2018real}. In this subsection, we will analyse the  accuracy of the proposed method in satisfying the final condition specified in Eq.~(\ref{FinalState}).

We firstly attempt to implement the landing error analysis by applying the proposed method to some pinpoint landings with initial conditions randomly generated in a defined range. However, we find that few landing cases succeed. This is probably caused by the fact the final condition specified in Eq.~(\ref{FinalState}) is fully constrained; in such case, small unfavorable trajectory dispersion could require the thrust vector to rapidly change, in order to satisfy the terminal position and speed constraints. Such rapid change may not be flyable for any lunar lander, resulting in the pinpoint problem having no solution \citep{lu2023propellant}. For this reason, $100$ initial conditions are randomly chosen in the dataset $\mathcal{D}$, and the proposed method is applied to guiding the lunar lander to the desired landing site. Any landing with a terminal speed error $V_f$ less than $5$ m/s is deemed as successful.

\begin{figure}[!htp]
\begin{center}
\includegraphics[scale=0.2]{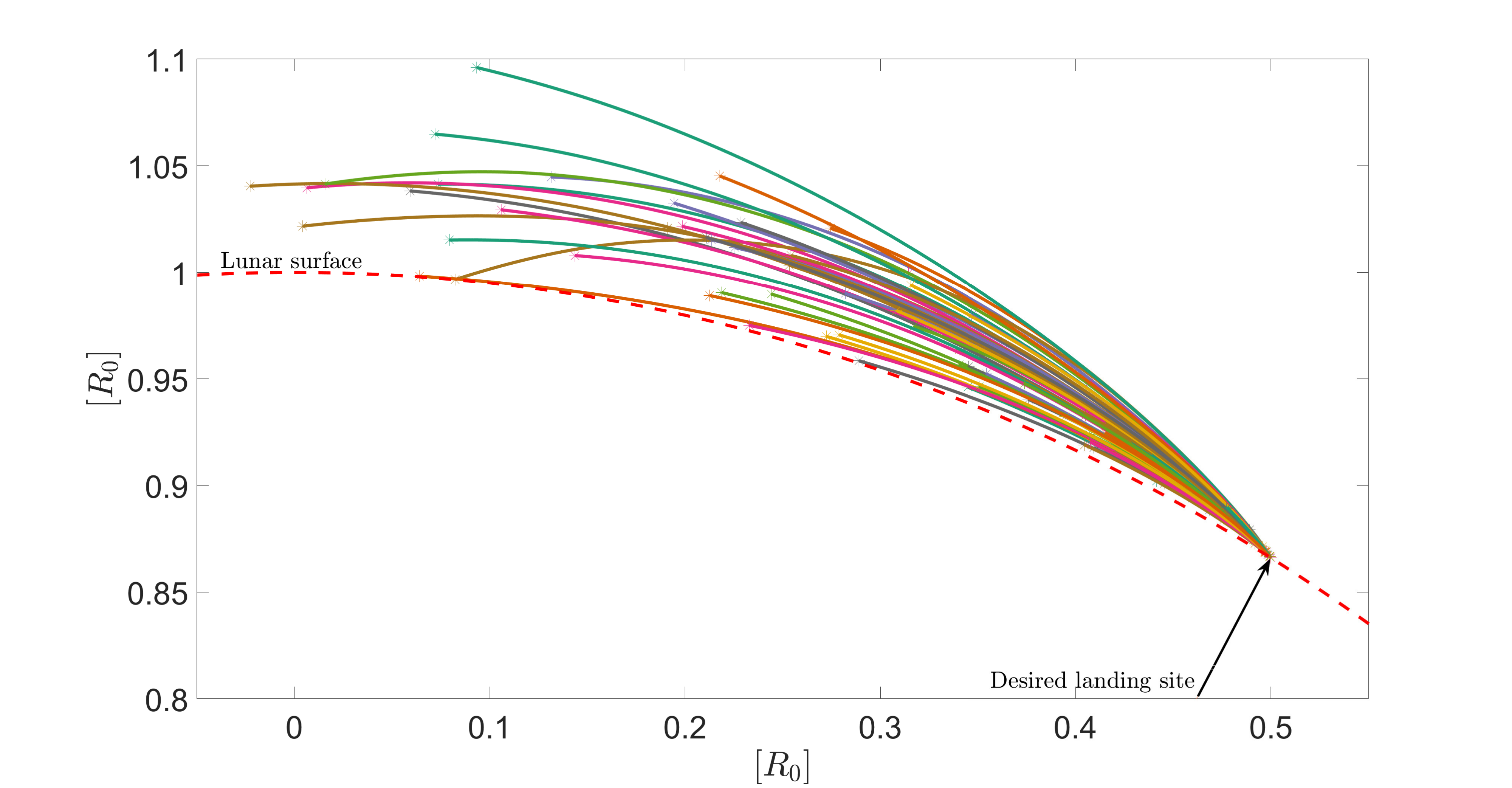}
\caption{NN-driven flight trajectories with different initial conditions randomly sampled in the dataset $\mathcal{D}$.}\label{Fig:100tra}
\end{center}
\end{figure}

\begin{figure}[!htp]
\centering
\begin{subfigure}[t]{7cm}
\centering
\includegraphics[width = 8cm]{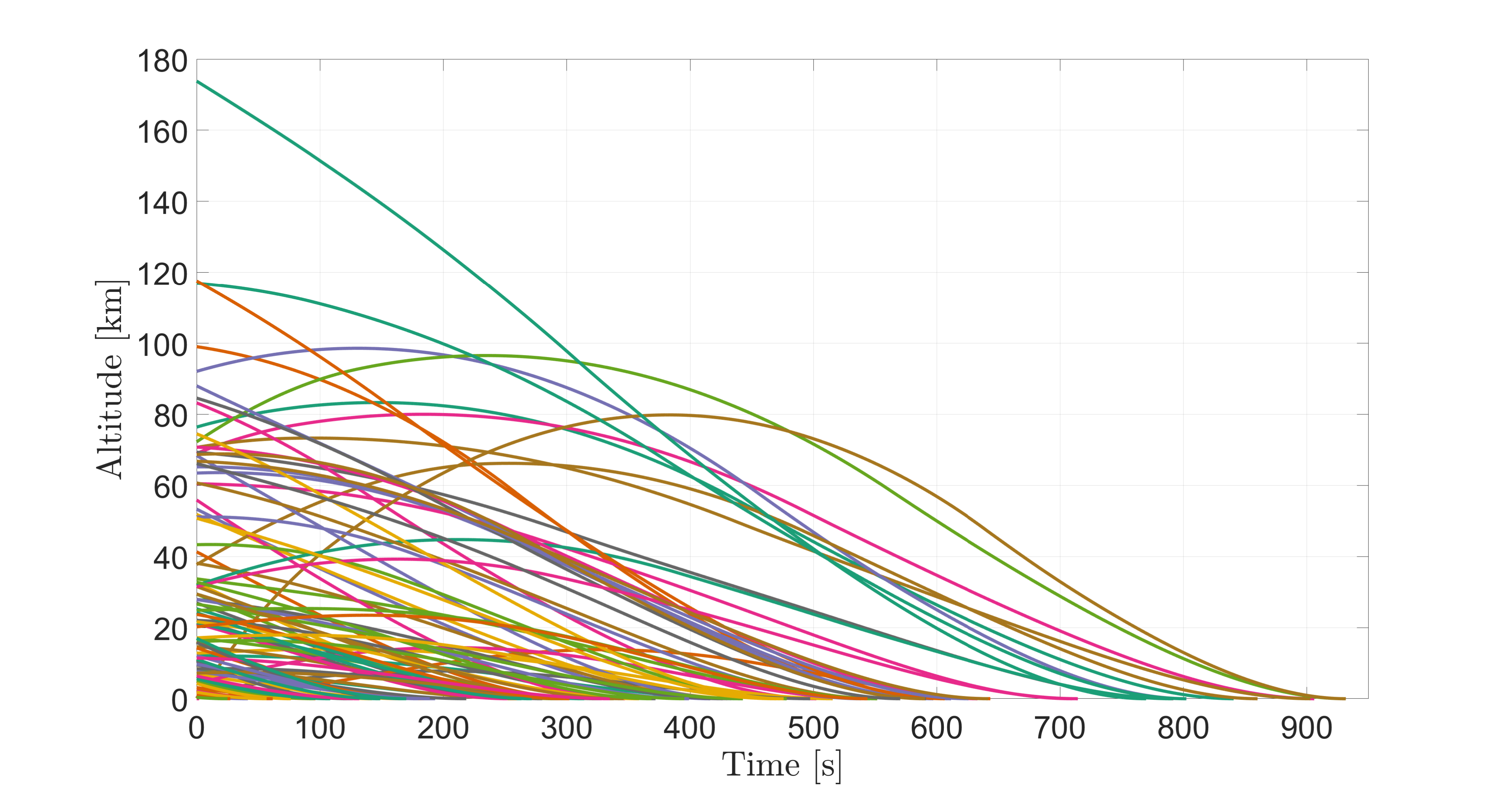}
\caption{Altitude profiles}
\label{Fig:100_alititude}
\end{subfigure}
~~~~~
\begin{subfigure}[t]{7cm}
\centering
\includegraphics[width = 8cm]{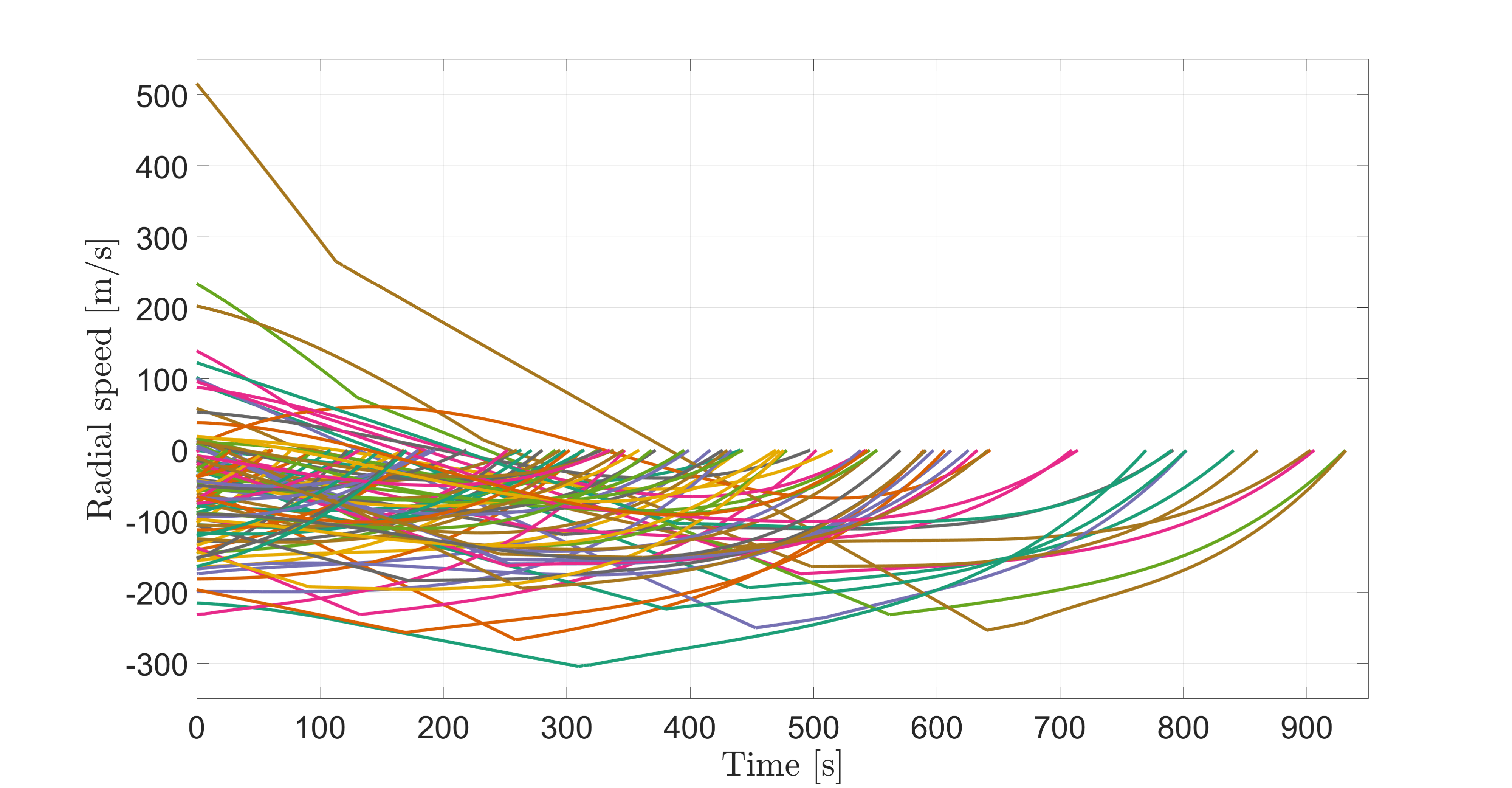}
\caption{Radial speed profiles}
\label{Fig:100_radialspeed}
\end{subfigure}
~~~~~
\begin{subfigure}[t]{7cm}
\centering
\includegraphics[width = 8cm]{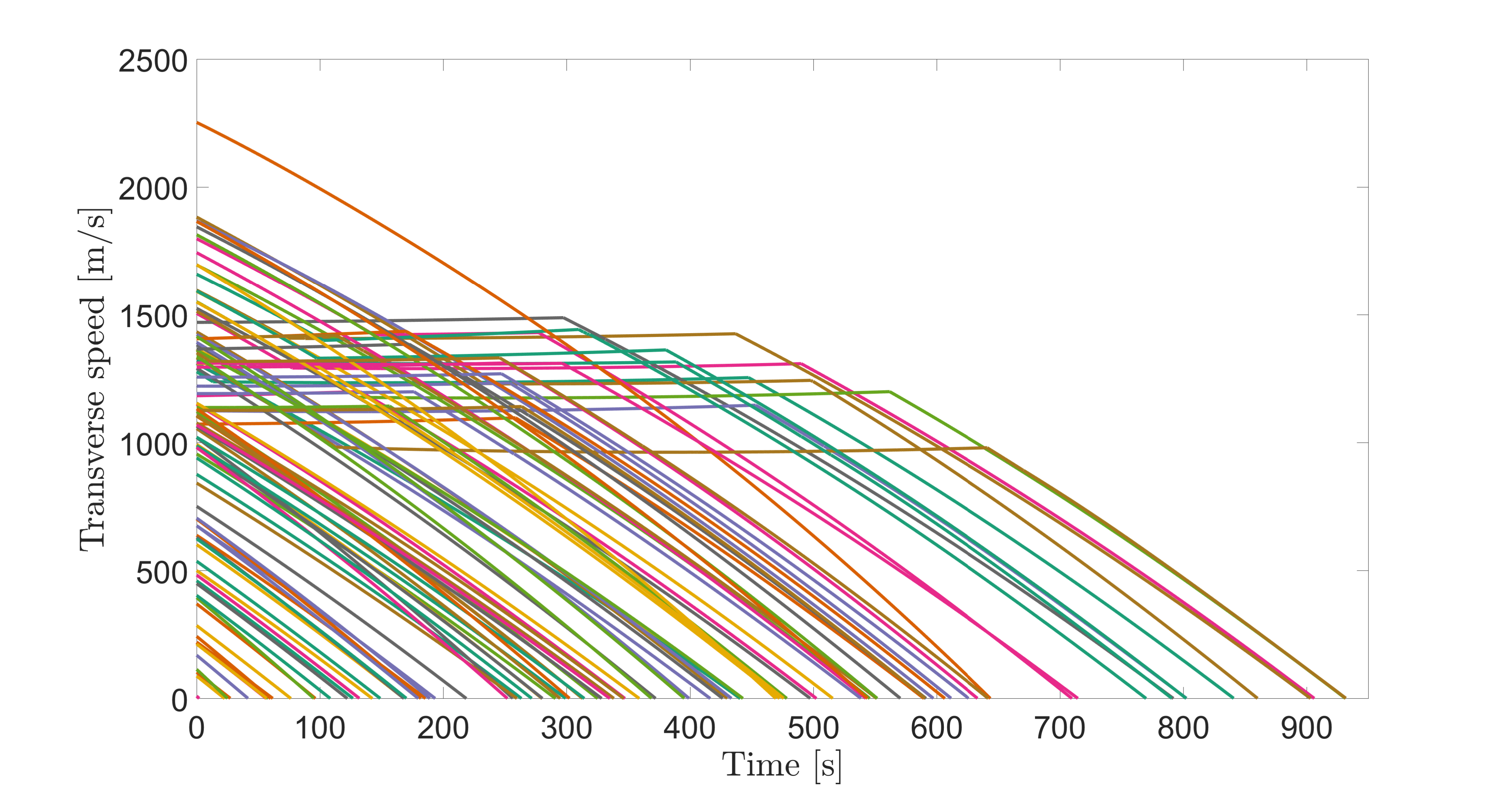}
\caption{Transverse speed profiles}
\label{Fig:100_tanspeed}
\end{subfigure}
~~~~~
\begin{subfigure}[t]{7cm}
\centering
\includegraphics[width = 8cm]{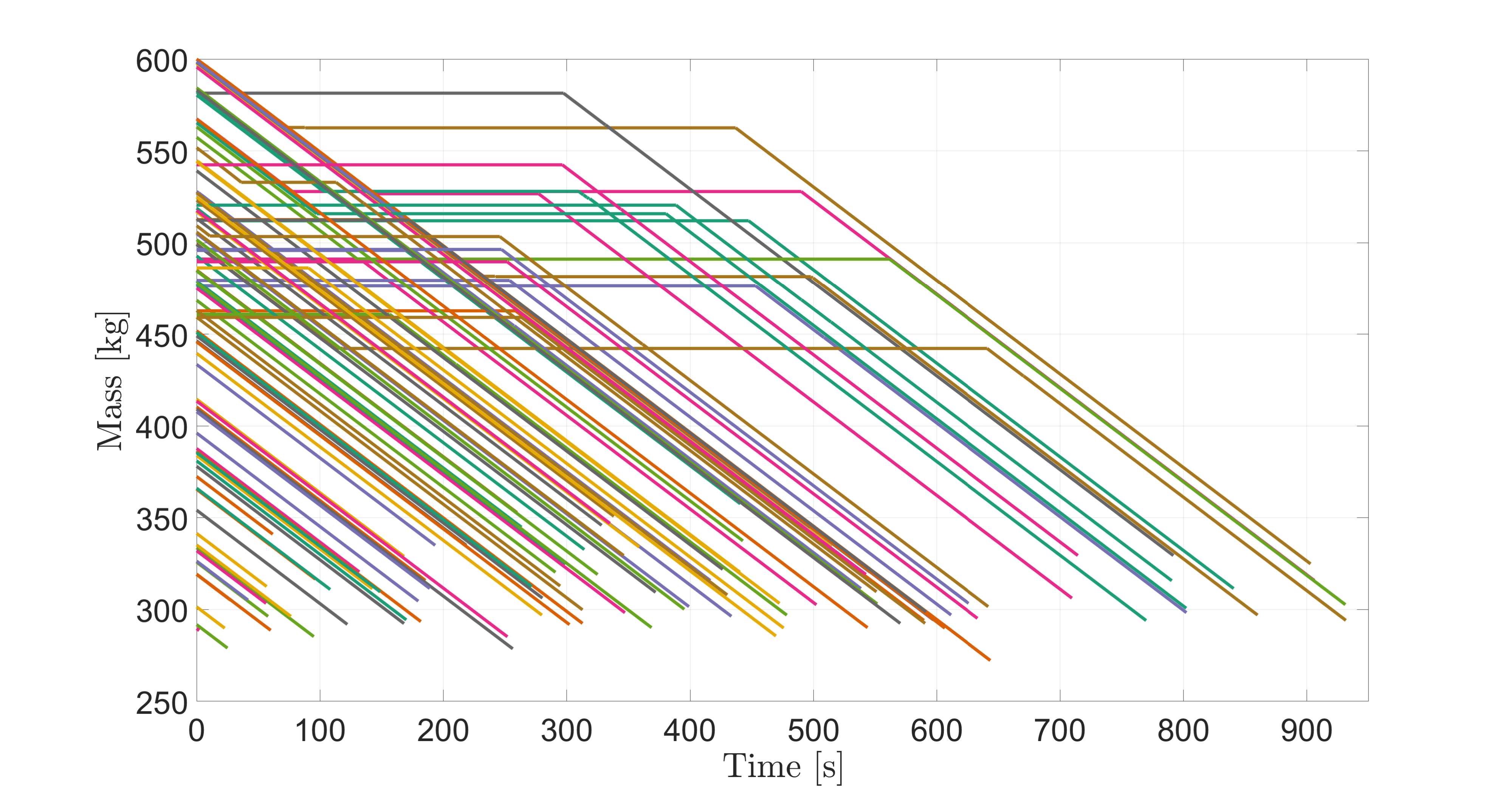}
\caption{Mass profiles}
\label{Fig:100_mass}
\end{subfigure}
\caption{Profiles of altitude, radial speed, transverse speed, and mass related to the NN-driven flight trajectories.}
\label{Fig:Profiles_100}
\end{figure}

Out of $100$ runs, a total of $100$ pinpoint landings are successful, and the corresponding flight trajectories are displayed in Fig.~\ref{Fig:100tra}, where the marker $*$ denotes the initial position of the lunar lander.
As it is shown in Fig.~\ref{Fig:100tra} that although the initial conditions of the lunar lander are quite different, the proposed method is still able to guide the lunar lander from the initial condition to the desired landing site. Fig.~\ref{Fig:Profiles_100} illustrates the profiles of altitude, radial speed, transverse speed, and mass corresponding to the $100$ successful cases. Notably, we can see from Fig.~\ref{Fig:100tra} that although the lunar lander is very close to the lunar surface at the beginning in some cases, the lunar lander is driven far away from the lunar surface at the initial stage because of high radial speed, as shown by the profiles of altitude and radial speed in 
Figs.~\ref{Fig:100_alititude} and \ref{Fig:100_radialspeed}, respectively; the landing phase is then followed by guiding the lunar lander to approach the desired landing site, in order to meet the final position and speed constraints. From Fig.~\ref{Fig:100_tanspeed}, it can be seen that the transverse speed of the lunar lander is generally decreasing during the landing phase; meanwhile, in some cases the transverse speed basically stays unchanged in the intermediate stage. This is because in such cases, the engine of the lunar lander is switched off and the impact of the gravitational force on the transverse speed is almost negligible compared to that of the engine thrust. In addition, we can observe that the engine is either kept on  or takes an ``on-off-on'' or ``off-on'' profile during the entire powered descent, as further indicated by the mass profiles in Fig.~\ref{Fig:100_mass}.

\begin{figure}[!htp]
\centering
\begin{subfigure}[t]{7cm}
\centering
\includegraphics[width = 8cm]{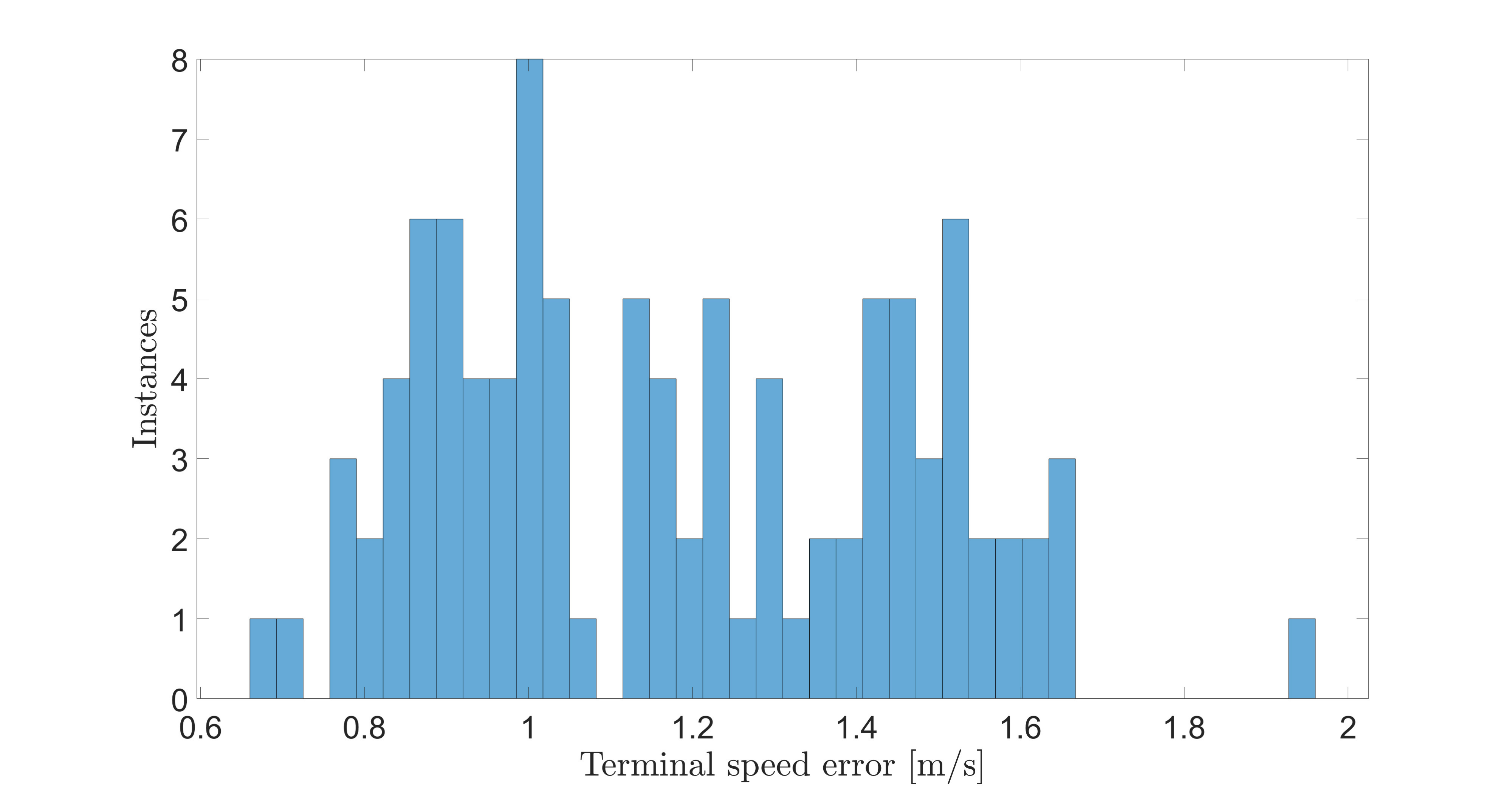}
\caption{Histogram of terminal speed error}
\label{Fig:100_terminalspeed}
\end{subfigure}
~~~~~
\begin{subfigure}[t]{7cm}
\centering
\includegraphics[width = 8cm]{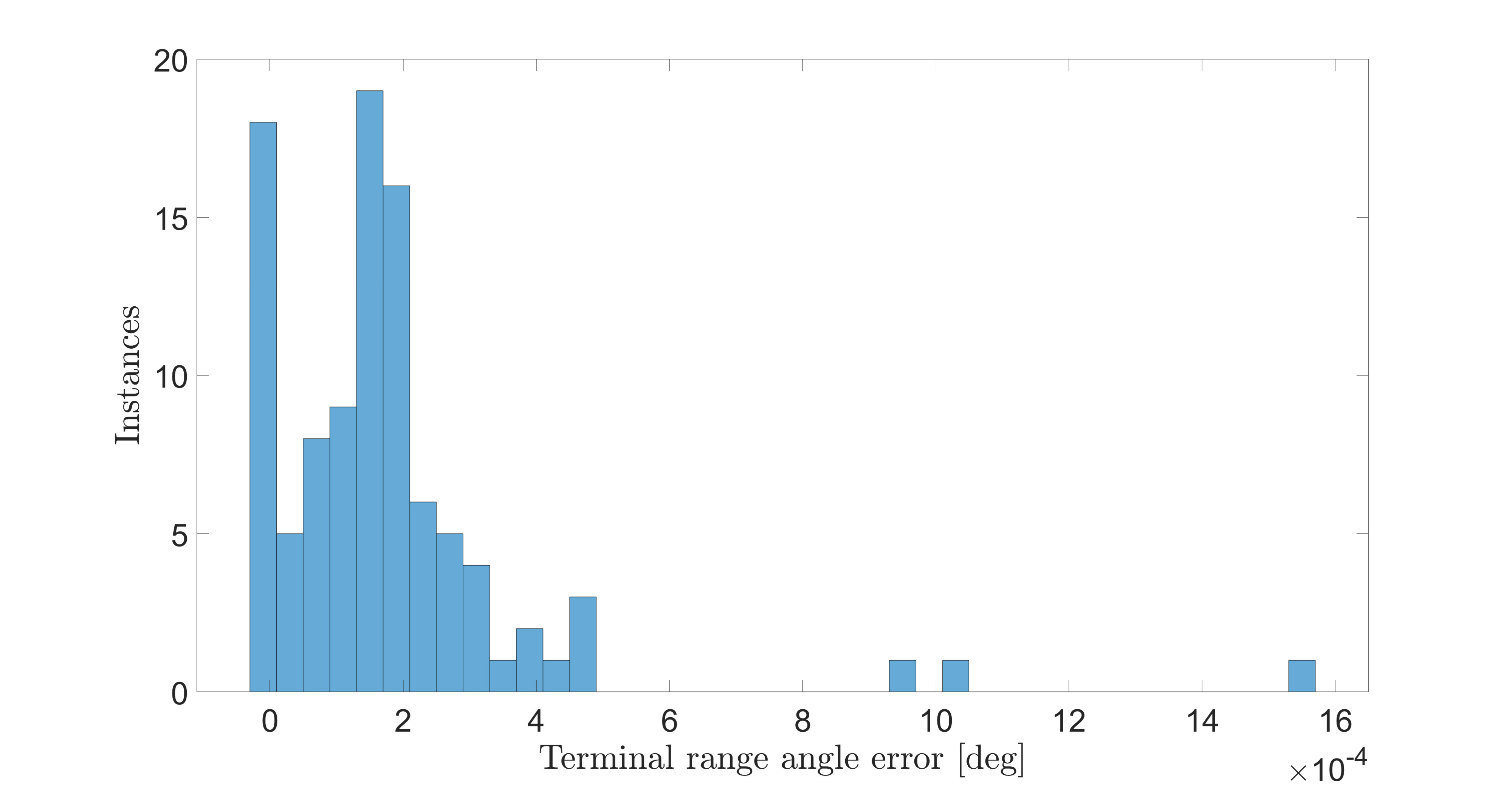}
\caption{Histogram of terminal range angle error}
\label{Fig:100_rangeangle}
\end{subfigure}
~~~~~
\begin{subfigure}[t]{7cm}
\centering
\includegraphics[width = 8cm]{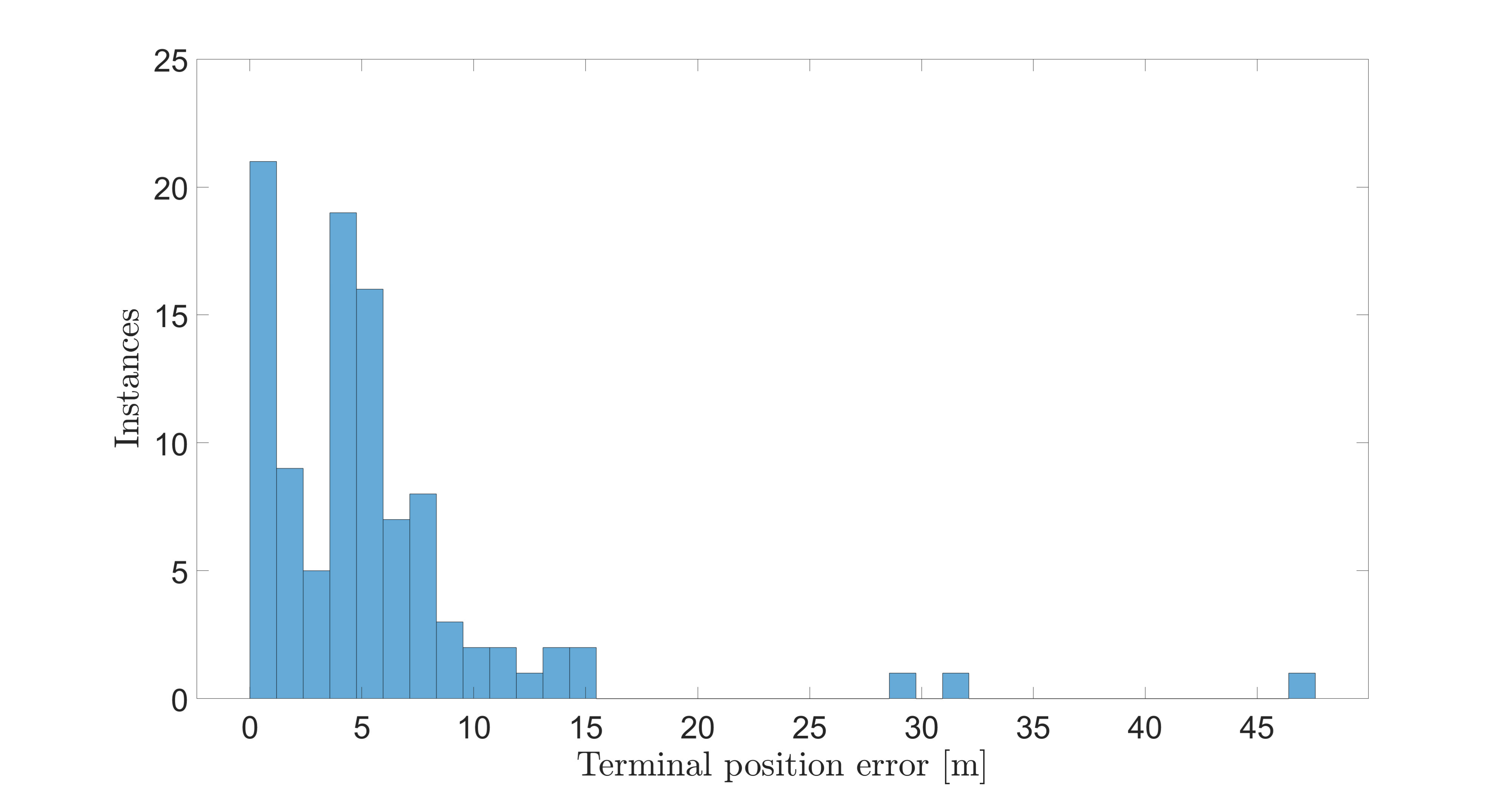}
\caption{Histogram of terminal position error}
\label{Fig:100_terminalposition}
\end{subfigure}
~~~~~
\begin{subfigure}[t]{7cm}
\centering
\includegraphics[width = 8cm]{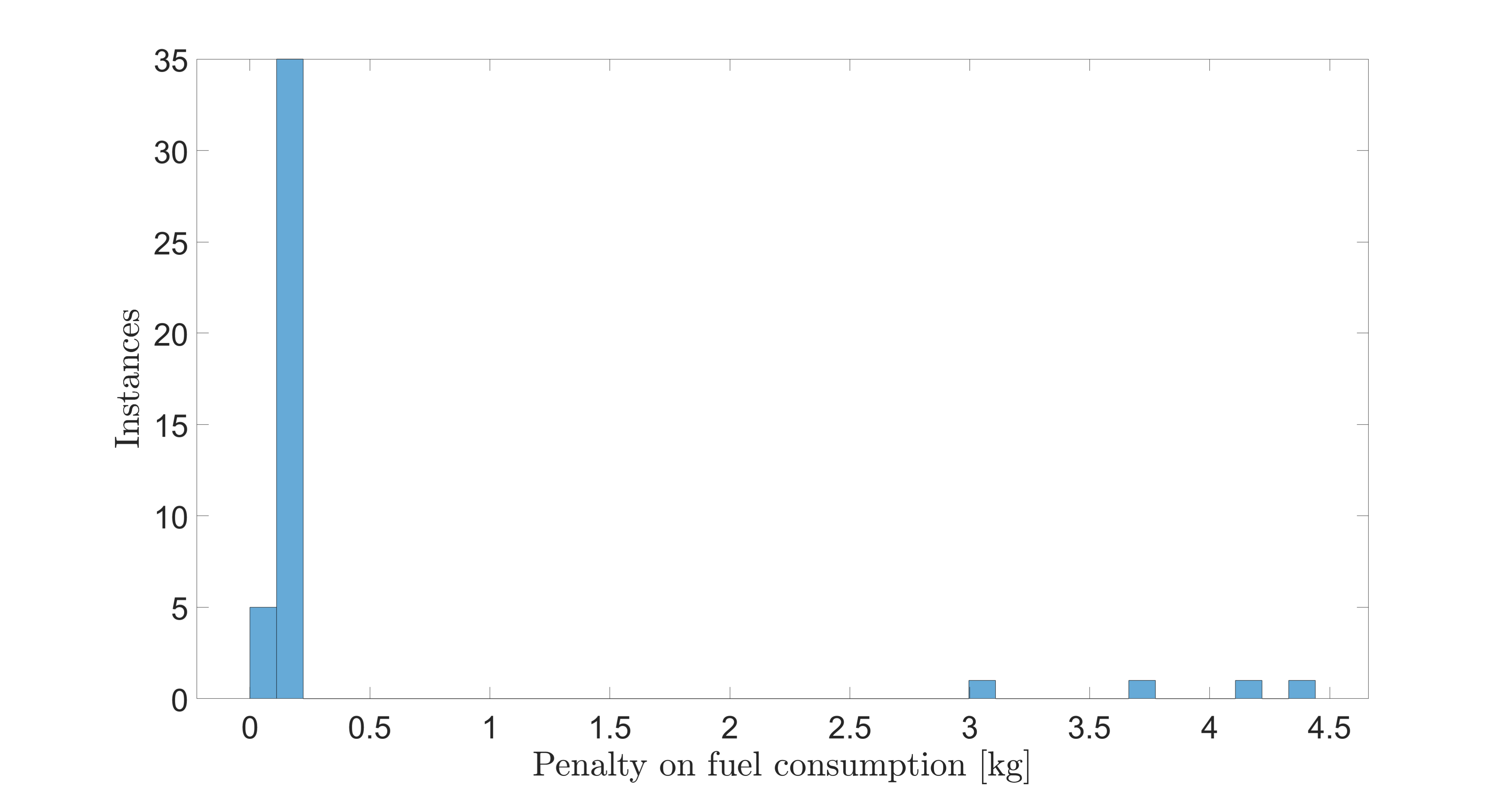}
\caption{Histogram of penalty on the fuel consumption}
\label{Fig:100_fuel}
\end{subfigure}
\caption{Error histograms of terminal speed, terminal range angle, terminal position, and fuel consumption penalty with 40 bins.}
\label{Fig:Error_histograms}
\end{figure}

At touchdown, the accuracy of the proposed method in satisfying the final condition specified in Eq.~(\ref{FinalState}) can be demonstrated by two key indicators, i.e., the terminal speed error $V_f$ and terminal range angle error $|\theta^{\mathcal{N}}_f|$ ($\theta^{\mathcal{N}}_f$ is the terminal range angle obtained from the proposed method).  The error histograms of these two indicators are depicted in Figs.~\ref{Fig:100_terminalspeed} and \ref{Fig:100_rangeangle}. It can be seen from Fig.~\ref{Fig:100_terminalspeed} that among the $100$ successful cases, the largest terminal speed error does not exceed $2$ m/s, and most of the successful cases have terminal speed errors concentrating within the range $[0.985,1.0175]$ m/s, indicating that the proposed method only results in a small terminal speed. Regarding Fig.~\ref{Fig:100_rangeangle}, the terminal range angle error does not exceed $0.0016$ deg. 
For clearer presentation, we define $e_p$ as the terminal distance between the landing site guided by the proposed method and the desired landing site, and  $e_p$ is computed by
\begin{align*}
e_p = \frac{2\pi R_0}{360}|\theta^{\mathcal{N}}_f|. 
\end{align*}
The error histogram of terminal position error $e_p$ is shown in Fig.~\ref{Fig:100_terminalposition}. Therefore, 
the proposed method is capable of guiding the lunar lander to the designated landing site with a small terminal position error. 

To assess the penalty on the fuel consumption caused by the proposed method, we apply the indirect method to solve the shooting function in Eq.~(\ref{EQ:TPBVP_law}) so as to find the optimal fuel consumption. Unfortunately, out of the $100$ landing cases, only $44$ cases converge for the indirect method, even if a homotopy on the smoothing constant in the optimal thrust magnitude is adopted. The histogram of penalty on the fuel consumption for these $44$ cases is presented in Fig.~\ref{Fig:Error_histograms}, from which we can see that even the largest penalty is less than $4.5$ kg, and the penalty on the fuel consumption for $40$ out of these $44$ cases is less than $0.25$ kg.
\section{Conclusions}\label{SE:conclusions}
Although Neural Networks (NNs) are widely applied in aerospace engineering to generate the optimal guidance command in real time, the conventional numerical optimization-based methods for acquiring the optimal state-guidance pairs are usually time consuming and may suffer from convergence issues. In addition, it has been shown in the literature that utilizing NNs to directly approximate a bang-bang control is quite challenging. To this end, we  parameterized the optimal trajectory based on the necessary conditions in virtue of Pontryagin's Minimum Principle, allowing for generating an optimal trajectory readily. Then, instead of directly approximating the bang-bang control via NNs, we attempted to approximate the continuous switching function. However, it was proven that the mapping from the flight state to the switching function is set-valued, at least in the final phase of the powered descent. The set-valued nature of such mapping poses an intricate challenge for training and using NNs. To eliminate the set-valued  mapping, a regularisation function was introduced. As a result, the approximation performance for the regualrised switching function was greatly enhanced. Meanwhile, another two well-trained NNs were able to predict the optimal thrust steering angle and time of flight given a flight state. Numerical simulations have shown that the proposed method could generate the fuel-optimal guidance in real time for the pinpoint landing with acceptable landing errors. Future research directions include generalization of the proposed method to pinpoint landing problem in three dimension and/or with more complex constraints, such as glide slope constraint and thrust pointing constraint.

\section*{Acknowledgments}
This research was supported by the National Natural Science Foundation of China under Grant No. 62088101.

\bibliography{refs}

\end{document}